\documentclass{amsart} 

\usepackage{tikz}
\usepackage[utf8]{inputenc} 
\RequirePackage{algorithm2e}
\SetAlFnt{\footnotesize}
\SetKwComment{Comment}{$\triangleright$\ }{}

\renewcommand\appendixname{Supplementary Material}
\usepackage{geometry} 
\usepackage{graphicx} 
\usepackage{hyperref}
\hypersetup{pdfauthor={Giacomo Aletti},pdftitle={Analytical confidence intervals for the number of different objects in data streams}}

\usepackage{amssymb,mathrsfs,bbm}
\usepackage{url}

\newcommand*\sref[1]{%
    S:\ref{#1}}

\newcommand*\seqref[1]{%
    (S:\ref{#1})}
    
\newcommand*\meqref[1]{%
    (M:\ref{#1})}

\newcommand{\hN}[1]{\mathbbm{h}_{#1}}

\theoremstyle{plain}
\newtheorem{theorem}{Theorem}[section]                                          
                          
\newtheorem{lemma}[theorem]{Lemma}

\theoremstyle{definition}
\newtheorem{definition}[theorem]{Definition}
\theoremstyle{remark}

\newtheorem*{esa}{Guiding example}

\usepackage{fancyhdr} 
\allowdisplaybreaks[4]

\title{Analytical confidence intervals for the number of different objects in data streams}
\author{Giacomo Aletti}
\address{G. Aletti,
ADAMSS Center, Universit\`a degli Studi di Milano, Milan, Italy}
\email{giacomo.aletti@unimi.it}

\begin{document}

\begin{abstract}
This paper develops a new mathematical-statistical approach to analyze
a class of Flajolet-Martin algorithms (FMa), and provides 
analytical confidence intervals for
the number $F_0$ of distinct elements in a stream, based on Chernoff bounds.
The class of FMa has reached a significant popularity in bigdata stream learning,
and
the attention of the literature has mainly been based on algorithmic aspects,
basically complexity optimality, while 
the statistical analysis of these class of algorithms has been
often faced heuristically.
The analysis provided here 
shows deep connections with mathematical special functions 
and with extreme value theory. 
The latter connection may help in explaining heuristic considerations, while
the first opens 
many numerical issues, faced at the end of the present paper. 
Finally, the algorithms are tested on an anonymized real data stream
and MonteCarlo simulations are provided to support our analytical choice in this context. 
\end{abstract}
\maketitle

\begin{footnotesize}
\noindent
\textbf{Acknowledgements.}
G. Aletti is a member of “Gruppo Nazionale per il Calcolo Scientifico (GNCS)” of the Italian “Istituto Nazionale di Alta Matematica (INdAM)”.

\noindent
\textbf{Competing interests.}
The author declares that he has no competing interests.

\noindent
\textbf{Availability of data and materials.}
All data, codes, and materials are available upon request.

\noindent
\textbf{Funding.}
This work has been partially supported by ADAMSS Center funds for Big Data research.
\end{footnotesize}

\section{Introduction}
Data streams \cite{stream00} are sequences of objects that cannot be
available for random access but must
be analyzed sequentially when they arrive and immediately discharged. 
Streaming algorithms process data streams and have reached a very rich audience since 
the last decades.
Typically, these kinds of algorithms have a limited time to complete
their processes and have access to limited amount of memory, usually logarithmic in the quantity of interest.

One of the main applications in streaming algorithms concerns the problem of counting 
\emph{the number $F_0$ of distinct elements} in a stream. 
Different solutions have been developed to estimate $F_0$ conserving memory space. 
\subsubsection*{State of the art.}
In \cite{Flajolet1985}, the authors develop the first algorithm for approximating $F_0$ 
based on hash functions. This algorithm was then formalized and made popular in \cite{Alon96}, where it was presented
the forefather of the class of algorithms that takes the name of \emph{Flajolet-Marin algorithms} (here, FMa). 
Three extensions in FMa were presented in \cite{FMext02}, together with a complete
description of the drawback and of the strength of the previous attempts. The first optimal (in complexity) algorithm 
has been proposed and proved in \cite{KanNelWoo10} and, nowadays, the FMa covers a lot of applications.
As only an example, in \cite{GanJea06}, an application with multiset framework is developed from
one of the most recent versions of FMa, 
and it estimates the number of  ``elephants'' in a stream of IP packets (see also \cite{Xi17}).
To summarize the state of the art, the typical sketch-based algorithms include 
PCSA \cite{Flajolet1985}, LinearCounting \cite{Whang90} (and MultiResBitmap as a generalization \cite{MulResBit06}),
MinCount \cite{FMext02}, LogLog \cite{DurFla03}, and HyperLogLog \cite{FlaFusGanMeu07} (see also
a recent generalization in \cite{Pe21}).

The FMa class of algorithms is essentially based on the following concept. 
When an object arrives from the stream, one (or more, independent) hash functions are applied to it, and then
the object is immediately discharged. 
The results of these hash functions are melted with what saved in memory (that has
a comparable size). The memory is updated, if necessary, with the result of this procedure, 
and then the process is ready for the next object.
The estimate of $F_0$ may be queried when necessary, and it is a function of the memory content only.

The key point is the fact that the central operation is made with a function which must be associative,
commutative and idempotent, so that multiple evaluations on the same object do not affect the final outcome,
which results in the combination of the hash values of the $F_0$ distinct objects.
A good candidate for such a function is the $\max$ function applied to 
a ``signature'' of each object, that is the core of such streaming algorithms. 
The same idea has recently used for other distributed algorithms (see \cite{Ale18} for simulation of
discrete random variables), 
where new entries or single changes should not make all the algorithm starts afresh.

\subsubsection*{Original Contribution}
As stated before, the main contributions in the study of FMa have concerned complexity problems, and
a deep mathematical-statistical approach has not yet developed, even if this class of algorithm is
probabilistic. 
This paper is a first attempt in this direction. The main contribution here is the analytical and numerical control
of FMa based on a pure mathematical statistic approach, while we
leave the measure of the goodness of the FMa to other studies 
(see \cite{ERT16} for a continuously updated work).
In particular, 
we give here analytical confidence intervals for the quantity of interest $F_0$. 
More precisely, 
we analyze an extension of the algorithms given above, and 
given the significance level ${\alpha}>0$, we will find $a,b > 0$, function of the memory content, such that
\begin{equation}\label{eq:chB1}
P( a \leq f(F_0) \leq b ) \geq \alpha,
\end{equation}
where $f$ is a given, strictly increasing, special function.
It is important to note that the approximations for $F_0$ as in \eqref{eq:chB1} given in literature are not satisfactory. 
In some situations, the asymptotic behavior of the interval is calculated through a Central Limit Theorem 
(see \cite{FlaFusGanMeu07}), 
but the huge skewness implicit in the algorithm variables (even in logarithmic scale) makes the Central Limit Theorem
questionable.
To overcome this observation,
Chebichev and Markov bounds are sometimes used to compute confidence intervals (see the papers cited in
\cite{KanNelWoo10}),
where the results are analyzed in terms of optimal complexity (in space and time) without exploiting
possible benefits in reducing the magnitude of the interval length.

These facts suggest us to not base the confidence interval on statistical asymptotic properties, but to build
analytical confidence intervals based on concentration inequalities.
In particular, we use here Chernoff bounds, and we give suitable approximations of the resulting inequalities. 
We show with MonteCarlo simulations that the analytical approximation does not affect the result
significantly. Moreover, we show that the same result derives from the use
of the Chernoff bounds on the limiting distribution that would be obtained with extreme value theory.

It is not surprising that some new analytical special functions appear in the analysis of the algorithm.
A particular modification of the analytical extension $\hN{1}(x)$ of the harmonic numbers function 
arises here as the mean value of a crucial quantity, so that $\hN{p}(\ln(2)F_0)$ is a quantity that appears in the paper. 
%

In addition, we discuss a numerical implementation of the analytical confidence intervals that can be run in real time. 
To do so, we analyze deeply all the relevant nonlinear problems that must be solved to build such confidence intervals. 
Then we provide the necessary numeric bounds to apply a new algorithm with a cubic rate of convergence, 
that has been tested successfully on a real anonymized data stream. 
As a byproduct, we give the algorithm that calculates the $\log$-shortest confidence interval
for $F_0$ based on the previous bounds. 

\subsubsection*{Organization of the paper.}
The paper is structured in the following way.
In the next Section~\ref{sec:description} we provide the quantities (parameters and statistics) used in the paper.
The description of both the streaming and the querying algorithms is given in the Section~\ref{sec:desc_alg}.
The main result, Theorem~\ref{thm:ris1}, is given
at the beginning of the Section~\ref{sec:apprY}, together with the connection
with the asymptotic results of the extreme value theory in Section~\ref{sec:extr_value}. 
The Section~\ref{sec:analyticalDiscussion} shows the goodness in the choice 
of the analytical approximations
given in the proof of the main theorem.
The algorithms given in this paper are tested on Twitter data 
and on a real anonymized 
data stream in Section~\ref{sec:realDataAnalysis}. 
In Section~\ref{sec:CompAsps} we face numerically some nonlinear equations related to the 
querying phases of the algorithm.
The mathematical properties of the special functions used in this paper,
the details of the 
proof of the main results, and the technicalities
needed to find lower and upper bounds 
contained in Section~\ref{sec:CompAsps} are left to
Supplementary Material \cite{Ale20SUPP}.
When necessary, the reference 
to the Supplementary Material are proceeded with a S, so that (S:A.$1$) will refer to the equation
\cite[(A.1)]{Ale20SUPP}.

\section{Description of the parameters and statistics of this paper}\label{sec:description}
The quantity $F_0$ denotes here the quantity of interest.
It gives the unknown number of distinct elements in
a real-time stream of possible repeating objects, and it is set as unknown parameter. 
The stream data is defined here as a sequence of objects $\{o_1, o_2, \ldots\}$.

We recall that FMa bases the $F_0$ estimate by
counting the maximum number of leading zeros in the hash values 
of the stream objects. One needs \(\log_2(F_0)+k\) bits in the hash function,
where the constant $k$ ensures a probability of the order of $\exp(-2^k)$ of having all bits
equal to $0$ in some hash values. 

In this paper the estimation is based on $c_0$ given independent hash functions
\(\{H_c,c=1,\ldots,c_0\}\). The main statistics of the first real-time phase
are extracted from the values that are resulting in applying
these functions on each object $o$ of the data stream.  
The results of the hash mapping
\(\{H_c(o),c=1,\ldots,c_0\}\) are used to
fill in-memory matrices $\mathbb{X} $
and $\mathbb{Z} $ of common size $2^{{r}_0}$ rows and $c_0$ columns (the total size
of such matrices will be denoted by $a_0 = 2^{{r}_0}c_0$).
The content of $\mathbb{X} $
and $\mathbb{Z} $ are then used during the querying phase to provide the confidence interval.
This memory data structure is a generalization of a 
HyperLogLog data structure (see \cite{DurFla03,KanNelWoo10,ERT16}).
The experimenter may choose the non-negative integer number $r_0$, together with 
another non-negative integer number $z_0$, 
to increases the accuracy of the estimates, at the cost to be sure that each hash function 
provides a sequence of bits 
longer than \(r_0+z_0+\log_2(F_0)+k\), with $k$ as above. 

\begin{esa}
In a word count streaming problem, the word \emph{pippo} is analyzed
and is mapped by first hash function $H_1$ to 
$H_1(\mathrm{pippo})= \mathrm{0xd012f681}$ (hexadecimal), that has a binary representation
given by 
\[
H_1(\mathrm{pippo})= {11010000000100101111011010000001}_2.
\] 
Then, with $r_0=4$ and $z_0= 6$, 
\begin{itemize}
\item
the first $r_0=4$ bits ${1101}_2=13$ of $H_1(\mathrm{pippo})$ are used
to build the first ``random'' number $R=1+13\in \{1,\ldots,2^{r_0}\}$, 
\item the successive $z_0=6$ bits ${000000}_2=0$ 
set the second quantity $Z=0\in \{0,\ldots,2^{z_0}\}$, 
\item the remaining bits $0100101\ldots$ 
are used to extract the number of the position of the first \emph{bit-one}: $X=2\in \{1,2,\ldots\}$.
\end{itemize}
The values of $R$, $X$ and $Z$ for the distinct objects of two datasets are plotted in Figure~\ref{fig:uniform}.
\end{esa}
Summing up, 
we denote by $H_c(o)$ the value of the $c$-th hash function applied to the object $o$, and 
it will consist
of a sequence of bits: $H_c(o)=(s_1,s_2, \ldots)$.
On this sequence, three statistics are extracted: $R=R(o,c)$ 
(from the first $r_0$ bits), $Z=Z(o,c)$ (from the subsequent
$z_0$ bits) and $X=X(o,c)$ (from the remaining bits). 
The quantities $X(o,c)$ and $Z(o,c)$ will update the elements $\mathbb{X}_{R,c}$ and
$\mathbb{Z}_{R,c}$, respectively, and then $R,X,Z$ are discharged.

During the second querying phase we build the confidence intervals. In this phase,
the central mathematical object are the statistics $\{Y_{r\,c},r=1,\ldots,2^{r_0},c=1,\ldots,c_0\}$. 
Each variable $Y_{r\,c}$ is a measurable function of the quantities $\mathbb{X}_{r,c}$ and
$\mathbb{Z}_{r,c}$, and the confidence interval at level $\alpha$
is made on the mean value $\mathcal{Y}$ of
these statistics.

\section{Description of the algorithm}\label{sec:desc_alg}
The \textbf{streaming algorithm} that updates $\mathbb{X}$ and $\mathbb{Z}$ in memory is given
in Algorithm~\ref{alg:storing_data}. 

\begin{algorithm}[t]
 \SetAlgoLined
\KwData{Data Stream of Objects $\{o_1,o_2, \ldots, \}$}
 \KwIn{$c_0$ hash functions, ${r}_0 \geq 0$ and $z_0 \geq 0$ small integers}
 \KwOut{Two matrices $\mathbb{X}$ and $\mathbb{Z}$ with $r_0= 2^{{r}_0}$ rows and $c_0$ columns}
 Set $\mathbb{X} \equiv 0$, $\mathbb{Z} \equiv 2^{z_0}-1 $ (binary)\;
\ForEach{$o$ in Stream}{
	\For{$c\leftarrow 1$ \KwTo $c_0$}{
		\tcc{compute the $c$-hash function on $o$, obtaining a sequence $(s_1,s_2,\ldots)$ of $0$ and $1$}
		$(s_1, s_2 , \ldots) \leftarrow H_c(o)$\;
		$R \leftarrow 1 + \sum_{r=1}^{{r}_0} s_r 2^{r-1}$ \Comment*[r]{$R \in \{1,\ldots,2^{{r}_0}\}$}
		$Z \leftarrow \sum_{z = 1}^{z_0} s_{{r}_0 + z} 2^{z_0-z} $ \Comment*[r]{$Z \in \{0,\ldots,2^{{z}_0}-1\}$}
		$X \leftarrow \inf\{n \geq 1 \colon s_{{r}_0 + z_0 + n} = 1\}$
		\Comment*[r]{{\footnotesize{$P(X +{r}_0 + z_0 > \text{length of hash}) \ll 1$}}}
		\uIf{$X > \mathbb{X}_{R\,c}$}{
			$\mathbb{X}_{R\,c} \leftarrow X$\;
			$\mathbb{Z}_{R\,c} \leftarrow Z$\;
		}\uElseIf{$X = \mathbb{X}_{R\,c}$}{
			$\mathbb{Z}_{R\,c} \leftarrow \min(Z,\mathbb{Z}_{R\,c})$\; 
		}
	}
	discharge $o$, $R$, $X$, $Z$\;
}
\caption{Streaming algorithm to store the data in memory. $\mathbb{X}$ is an integer-valued matrix,
whose values are of the order of $\log_2(F_0)$, while $\mathbb{Z}$ has 
values in $0, \ldots, 2^{z_0}-1$}\label{alg:storing_data}
\end{algorithm}
The flow of information is as follows. An object $o$ arrives in the stream data.
Each hash function $H_c$ applied to $o$ produces a sequence $(s_1,s_2,\ldots)$ of bits, from which
we extract 
$R = 1 + \sum_{r=1}^{{r}_0} s_r 2^{r-1}$, $Z = \sum_{z = 1}^{z_0}
s_{{r}_0 + z} 2^{z_0-z} $
and $X = \inf\{ n \geq 1 \colon s_{{r}_0 + z_0 + n} = 1\}$:
\begin{equation}\label{eq:Hco}
H_c(o) = 
\underbrace{0 1 \cdots 10 1}_{R \in \{1,\ldots,2^{{r}_0}\}}^{{r}_0\text{ bits}}
\underbrace{1 0 \cdots 0 1}_{Z}^{z_0\text{ bits}}
\underbrace{0 0 \cdots 00 0 1}_{X \in \{1,2,\ldots\}}^{X\text{ bits}}
\underbrace{0 1 1 0 1 0 0 0 \cdots}_{\text{not used}}
\end{equation}
The data are then updated according to the following procedure:
\begin{description}
\item[if $X<\mathbb{X} _{R\,c}$] do nothing;
\item[if $X>\mathbb{X} _{R\,c}$] set $\mathbb{X}_{R\,c} =  X$ and $\mathbb{Z}_{R\,c} = Z$;
\item[if $X=\mathbb{X} _{R\,c}$] 
set $\mathbb{Z} _{R\,c} = \min(\mathbb{Z} _{R\,c} , Z)$.
\end{description}

\begin{esa}[Continued]
With the guiding example started in the previous section, the result of the $c=1$-st hash function applied to the word \emph{pippo}
($R=14$, $Z=0$ and $X=2$) will cause a comparison with the content of $\mathbb{X} _{R=14,c=1}$ and $\mathbb{Z} _{R=14,c=1}$, and then
\begin{description}
\item[if $2<\mathbb{X} _{14,1}$] do nothing;
\item[if $2>\mathbb{X} _{14,1}$] set $\mathbb{X}_{14,1} = 2$ and $\mathbb{Z}_{14,1} = 0$;
\item[if $2=\mathbb{X} _{14,1}$] 
set $\mathbb{Z} _{14,1} = \min(\mathbb{Z} _{14,1} , 0)$.
\end{description}
\end{esa}

The \textbf{querying algorithm} first produces the matrix $\mathbb{Y}=\{{Y}_{r\,c},r= 1, \ldots , 2^{{r}_0}, c=1, \ldots, c_0\}$ 
with the contents of $\mathbb{X}$ and $\mathbb{Z}$: 
\begin{equation}\label{eq:fromXZtoY}
{Y}_{r\,c} = \mathbb{X}_{r\,c} - \log_2 (1+ 2^{-z_0}\mathbb{Z}_{r\,c} ),
\end{equation}
see
Algorithm~\ref{alg:querying_data}. 
Then the arithmetic mean $\mathcal{Y}$ of the $a_0 = c_0 2^{{r}_0}$ entries of
$\mathbb{Y}$ is evaluated to build a $\alpha$ confidence interval.
\begin{algorithm}[t]
 \SetAlgoLined
 \KwIn{$\mathbb{X}$ and $\mathbb{Z} $, output of Algorithm~\ref{alg:storing_data}}
 \KwOut{$\mathbb{Y} = \{ {Y}_{r\,c}, r= 1, \ldots , 2^{{r}_0}, c=1, \ldots, c_0\} $}
 Set $\tilde{Y} = 0$\;
\For{$c\leftarrow 1$ \KwTo $c_0$}{
	 \For{$r\leftarrow 1$ \KwTo $2^{{r}_0}$}{
		$y \leftarrow 2^{-z_0} \mathbb{Z}_{r\,c} $
		\Comment*[r]{$y \in [0, 1-2^{-z_0}] \ \Rightarrow \ (1+y) \in [1,2)$}
		${Y}_{r\,c} \leftarrow \mathbb{X}_{r\,c}  - \log_2(1+y)$
		\Comment*[r]{{\footnotesize{$\mathbb{X}_{r\,c}  - \log_2(1+y) \in (\mathbb{X}_{r\,c}-1,\mathbb{X}_{r\,c}]$}}}
	}
}
return $\mathbb{Y} =  ({Y}_{r\,c})_{r= 1, \ldots , 2^{{r}_0}, c=1, \ldots, c_0}$\;
\caption{Querying algorithm to extract $\mathbb{Y}$, starting from
the memory content $\mathbb{X}$ and $\mathbb{Z}$ given in
Algorithm~\ref{alg:storing_data}}\label{alg:querying_data}
\end{algorithm}
As an example, in Algorithm~\ref{alg:conf_int}, we compute
a $\alpha$-confidence interval for $F_0$
of the form $(0,\mathrm{upper})$, 
based on the Theorem~\ref{thm:ris1}. 
\begin{esa}[Continued]
Again, if we use the guiding example and we suppose that $\mathbb{X}_{14,1} = 2$ and $\mathbb{Z}_{14,1} = 0$,
we obtain the quantity ${Y}_{14,1} = 2 -\log_2 (1+2^{-6}\cdot 0) = 2$. Note that we always have that 
$1 \leq 1+2^{-z_0} Z_{r\,c} < 2$ which implies that $\mathbb{X}_{r\,c}-1 < {Y}_{r\,c} \leq \mathbb{X}_{r\,c}$.
The values of ${Y}_{r\,c}$ for two datasets are plotted in Figure~\ref{fig:uniform} (bottom-right).
\end{esa}
\begin{algorithm}[bht]
 \SetAlgoLined
 \KwIn{1) $\mathbb{Y} = \{ {Y}_{r\,c}, r= 1, \ldots , 2^{{r}_0}, c=1, \ldots, c_0\} $, 
 output of Algorithm~\ref{alg:querying_data}.\\2) the 
 confidence $\alpha\in (0,1)$ -usually $\alpha\in [0.9,0.995]$-}
 \KwOut{A $\alpha$ confidence interval for $F_0$ of the form $(0, \mathrm{upper})$}
Set $y = -\log(1-\alpha)/(2^{{r}_0}c_0)$\;
Set ${x} \leftarrow InvAlphaMinus(y)$ \tcc*{Solve (in $x$) the problem 
$y - ((x-\gamma)t_- - \ln(\Gamma(1+t_-))) =0$, with $\psi(1+t_-)= x-\gamma $}
Set $\hat{y} \leftarrow 0$\;
\For{$c\leftarrow 1$ \KwTo $c_0$}{
	 \For{$r\leftarrow 1$ \KwTo $2^{{r}_0}$}{
		$\hat{y} \leftarrow \hat{y} + {Y}_{r\,c}$.
	}
}
$\mathcal{Y} \leftarrow \hat{y} /(2^{{r}_0}c_0) $\;
Set $z \leftarrow \mathcal{Y} \log(2) + x + 2^{-z_0} $\;
Set $p_0 \leftarrow 2^{-{r}_0}$\;
return $\mathrm{upper} =  invHpM (z,p_0) $
 \Comment*{Solve (in $x$) the problem  $z - \hN{p_0} (x)  =0$}
\caption{Querying algorithm that builds a $\alpha$-confidence interval for $F_0$
of the form $(0,\mathrm{upper})$, 
based on the Theorem~\ref{thm:ris1}}\label{alg:conf_int}
\end{algorithm}

Finally, note that the data structure becomes that of \cite{ERT16}  when $c_0=1$ and $z_0=0$ 
(the content of $\mathbb{Z}$ is not significant and the update reduces to
$\mathbb{X}_{R\,c} \leftarrow \max(X,\mathbb{X}_{R\,c})$, without the if-else loop).
When, in addition, 
${{r}_0}=0$ the data structure reduces to the original
one \cite{FlaMar85}.

\subsection{Mathematical and Statistical analysis of the algorithm}\label{sec:stat_descr}
Given any object $o$ in the data stream, the streaming algorithm given in Algorithm~\ref{alg:storing_data}
extracts three measurable statistics $R$, $X$ and $Z$.
The first one is used to augment artifically the number of recorded statistics as in
\cite{ERT16}, while the latter ones deserve a more accurate explanation.
Take two objects $o_1$ and $o_2$, and assume that we collect $(R_1,X_1,Z_1)$ from
the first object, $(R_2,X_2,Z_2)$ from
the second one with the $c$-th hash function. 
If, by chance, $R_1=R_2=r$, then the contribution of these two objects to $\mathbb{Y}$
in the subsequent Algorithm~\ref{alg:querying_data} will be
\[
Y_{r\,c} = \max \big( {X}_{1} - \log_2 (1+ 2^{-z_0} {Z}_1 ) , {X}_{2} - \log_2 (1+ 2^{-z_0} {Z}_2 ) \big)
\]
as a consequence of \eqref{eq:fromXZtoY} and of the definition of $\mathbb{X}$ and
$\mathbb{Z}$. The $\max$ function here is the core of this algorithm, 
being a binary operation 
that has associativity, commutativity, and idempotence properity.
Algebraically speaking, 
a set $S$ with such a binary operation $\circ$ is called \emph{semilattice}. The key point is that semilattices $(S,\circ)$ 
are one-to-one
related to partially ordered relations $(S,\geq)$:
\(a \geq b \iff a \circ b = a\), so that they induce set operation instead point ones. In other simpler words,
when you evaluate the semilattices operator on different, even
repeated objects, the result is independent of the order and of the repetitions of the objects (as the $\max$ function does).
This fact is a mathematical key point when you want to estimate a function of the different objects without registering
the different objects you have seen so far. As a direct consequence,
%
the \emph{Algorithm~\ref{alg:storing_data} may be thought as applied only once  to 
each of the $F_0$ different objects}. 

\medskip

From a statistical point of view,
we will assume that each hash function generates an independent sequence of  bits that are equally distributed among all
the possible outcomes. In other words, we assume that the set $\{H_c(o), c=1,\ldots,c_0, o \text{ different objects}\}$
is made by a sequence of independent and identically distributed vectors of bits, each vector having 
bit components independent and equally distributed on $\{0,1\}$.  The sequence of bits $s_i$ in \eqref{eq:Hco}
is hence distributed as a Bernoulli of parameter $1/2$, and it is independent from the others.
%
%
%
%
%
%
%

Summing up, for each hash function $H_c$ and any object $o$ belonging to data stream,
the three statistics $R=R(c,o)$, $X=X(c,o)$ and $Z=Z(c,o)$ are collected, and the matrices $\mathbb{X}$
and $\mathbb{Z}$ updated. Then, during the querying phase, the statistics
\begin{equation}\label{eq:YrcMAT}
Y_{r\,c} = \max_{o\colon R(c,o)=r} \big( X(c,o) - \log_2 (1+ 2^{-z_0} Z(c,o) ) \big)
\end{equation}
is computed.

We now recall that, by definition, 
$2^{-z_0} Z(c,o) = 2^{-z_0} \sum_{z = 1}^{z_0} s_{{r}_0 + z} 2^{z_0-z}= \sum_{z = 1}^{z_0} s_{{r}_0 + z} 2^{-z} $.
This quantity may be seen as a truncated series. 
We complete the bit sequence $(s_{{r}_0 + 1},\ldots s_{{r}_0 + z_0}) $ and we form 
an i.i.d.\ sequence of equally distributed 
bits $(s^*_{1},\ldots s^*_{z_0}, s^*_{z_0+1} \ldots)$, 
where $s^*_{z} = s_{{r}_0 + z} $ if $z\leq z_0$. With this notation 
\[
2^{-z_0} Z(c,o) = \sum_{z = 1}^{z_0} s^*_{z} 2^{-z} ,
\]
the random variable
\[
\bar{Z}(o,c) = \sum_{z = 1}^{\infty} s^*_{z} 2^{-z}
\]
is uniformly distributed on $(0,1)$ and $0 \leq \bar{Z}(o,c) - 2^{-z_0} {Z}(o,c) < 2^{-z_0} $. 
More remarkable, if we denote by
\[
\bar{Y}(o,c) = \big( X(c,o) - \log_2 (1+ \bar{Z}(o,c) ),
\]
then
the random variable
\begin{align*}
\bar{U}(o,c) & =  2^{-\bar{Y}(o,c)} 
 = 2^{-X(c,o)} \big( 1 + \sum_{z = 1}^{\infty} s^*_{z}  2^{-z} \big) 
  = 2^{-X(c,o)} + 2^{-X(c,o)}\sum_{z=1}^{\infty} s^*_{z}  2^{-z}  \\
 & 
 = \sum_{x=1}^{X(c,o)} s_{{r}_0+{z}_0+x}2^{-x} + 
 \sum_{z=1}^{\infty} s^*_{z}   2^{-z+X(c,o)} ,
\end{align*}
is uniformly distributed on $(0,1)$, which immediately implies that
$\bar{Y}(o,c) = -\log_2(\bar{U}(o,c)) = -\frac{\log(\bar{U}(o,c))}{\log(2)}$
is an exponential random variable with parameter $\lambda_0 = \log (2)$.
The fact here is that, instead of measuring 
$\bar{Y} (o,c) 
$, 
we can only collect 
$X(c,o) - \log_2 (1+ 2^{-z_0} Z(c,o) )$, due to computational limitations,
and this introduces a further bias.
If we could have measured $\bar{Y} (o,c) $, the quantity \eqref{eq:YrcMAT} would have been
\begin{align*}
\bar{Y}_{r\,c} & = \max_{o\colon R(c,o)=r} \big( X(c,o) - \log_2 (1+ \bar{Z}(c,o) ) \big)
= \max_{o\colon R(c,o)=r} \big( \bar{Y}(o,c) \big)
\end{align*}
that is not too far from ${Y}_{r\,c}$, since we always have that 
$0< {Y}_{r\,c}  - \bar{Y}_{r\,c} < \tfrac{2^{-z_0}}{\lambda_0}$ (see \cite[Section~\sref{sec:Ybounds}]{Ale20SUPP}).
Finally, since
\[
\bar{Y}_{r\,c}
= \max_{o\colon R(c,o)=r} \big( \bar{Y}(o,c) \big)
 = 
\max_{\substack{o\colon R(c,o)=r\\ o \text{ different objects}}} \big( \bar{Y}(o,c) \big) ,
\]
the independence of the hash functions and of their results on different objects implies that
$\{\bar{Y}_{r\,c},r=1, \ldots,2^{-r_0},c=1, \ldots,c_0\}$ are a collection of independent random variables,
each of one being distributed as the maximum of a random number $m_{r\,c}$ of independent exponential random variables,
where
\[
m_{r \, c} = \#\big\{ o \in \{\text{$F_0$ different objects} \}\colon R(o,c) = r\big\}.
\]
It is obvious that, for any fixed $c$, $\sum_{r=1}^{2^{r_0}} m_{r \, c} = F_0$ and, moreover,
since $R= R(o,c) $ is uniformly distributed on $1,\ldots,2^{r_0}$, then 
the $c_0$ random vectors $\{\boldsymbol{m}_c =  (m_{1\, c}, \ldots, m_{2^{{r}_0}\, c}), c= 1,\ldots,c_0\}$ are
distributed as multinomial vectors of parameters $ F_0$ and $2^{-{r}_0}$, and independent of each other. 
%

We have proved the following result.
%
\begin{lemma}\label{cor:Y_rc}
There exists a family
\[
\big\{\bar{Y}(o,c), o \in \{\text{$F_0$ different objects}\}, 
c \in \{1, \ldots, c_0\} 
\big\}
\] 
of independent and identically distributed random variables
with exponential distribution of parameter 
$\lambda_0 = \log 2$, such that, if we define,
\[
\bar{Y}_{r\,c} = \max_{\{o \colon R(o,c)= r\}} ( \bar{Y}(o,c) ),
\]
then, uniformly in $r$ and $c$, 
\[
0 < {Y}_{r\,c}  - \bar{Y}_{r\,c} \leq \tfrac{2^{-z_0}}{\lambda_0},
\]
where each ${Y}_{r\,c}$ is defined in \eqref{eq:YrcMAT}. Moreover, 
for any fixed $c\in \{1, \ldots, c_0\}$, define
\[
m_{r \, c} = \#\big\{ o \in \{\text{$F_0$ different objects} \}\colon R(o,c) = r\big\}.
\]
Then
the random vectors $\{\boldsymbol{m}_c =  (m_{1\, c}, \ldots, m_{2^{{r}_0}\, c}), c= 1,\ldots,c_0\}$ are
i.i.d, distributed as multinomial vectors of parameters $ F_0$ and $2^{-{r}_0}$. 
Conditioned on $\boldsymbol{m}_c $, the random variables 
$\{\bar{Y}_{r\,c} , r= 1 , \ldots, 2^{r_0}\}$
are independent. 
\end{lemma}

\section{Confidence interval for $F_0$}\label{sec:apprY}
The main result of this section is the construction of a analytic confidence interval for $F_0$,
based on $\mathbb{Y}$ explained in the previous section. This interval is based on some 
special functions. The interested reader may find details in \cite[Section~\ref{app:Spec_func}]{Ale20SUPP}.
\begin{theorem}\label{thm:ris1}
Let 
$\mathbb{Y}$ 
be collected as in Section~\ref{sec:desc_alg}, and define
\[
\mathcal{Y} = \frac{\sum_{r=1}^{2^{{r}_0}}\sum_{c=1}^{c_0} Y_{r\, c}}{ 2^{{r}_0} c_0}.
\]
Then 
\begin{align*}
&\big(\hN{p_0}^{-1} (\lambda_0\mathcal{Y} - h_d)  , +\infty\big)
\\
&\big( 0 , \hN{p_0}^{-1} (\lambda_0\mathcal{Y} + h_u+\tfrac{2^{-z_0}}{\lambda_0})  \big)
 \\
& \big( \hN{p_0}^{-1}  (\lambda_0\mathcal{Y}  - h_d)  , \hN{p_0}^{-1}  (\lambda_0\mathcal{Y} + h_u+\tfrac{2^{-z_0}}{\lambda_0})  \big)
\end{align*}
are confidence intervals for the unknown parameter $F_0$, where
\begin{itemize}
\item $p_0=2^{-{r}_0}, \lambda_0 = \log(2)$; 
\item the function $\hN{p_0}:\mathbb{R}_+\to\mathbb{R}_+$ is defined as
\[
\hN{p}(x) = \int_0^1 \frac{1-(1- p + pt)^x}{1-t} \,dt , 
\]
\item the levels of confidence are ${\alpha}_+$, ${\alpha}_-$, and $({\alpha}_++{\alpha}_-)$ respectively, where
\begin{align*}
{\alpha}_+ & = 1 - \exp\Big( -2^{{r}_0}c_0 \big[ (h_d+\gamma) t_+ - \ln \Gamma(1-t_+) \big]\Big), \qquad t_+ = 1-\psi^{-1}(-h_d-\gamma);
\\
{\alpha}_- & = 1 - \exp\Big( -2^{{r}_0}c_0 \big[  (h_u-\gamma)  t_- - \ln \Gamma(1+t_-) \big]\Big), \qquad t_- = \psi^{-1}(h_u-\gamma)-1;
\end{align*}
$\gamma$ is the Euler constant and $\psi$ is the digamma function. 
\end{itemize}
\end{theorem}
\begin{proof}[Sketch of the proof of Theorem~\ref{thm:ris1}]
We first note that, by Lemma~\ref{cor:Y_rc}, if we define
\begin{equation}\label{eq:barY}
\bar{\mathcal{Y}} = 
\frac{\sum_{r=1}^{2^{{r}_0}}\sum_{c=1}^{c_0} \bar{Y}_{r\, c}}{ 2^{{r}_0} c_0},
\end{equation}
then $ 0 \leq \mathcal{Y} - \bar{\mathcal{Y}}  < \tfrac{2^{-z_0}}{\lambda_0}$, and then it is sufficient to prove that
\begin{align*}
&\big(\hN{p_0}^{-1} (\lambda_0\bar{\mathcal{Y}} - h_d)  , +\infty\big)
\\
&\big( 0 , \hN{p_0}^{-1} (\lambda_0\bar{\mathcal{Y}} + h_u)  \big)
 \\
& \big( \hN{p_0}^{-1}  (\lambda_0\bar{\mathcal{Y}} - h_d)  , \hN{p_0}^{-1}  (\lambda_0\bar{\mathcal{Y}} + h_u)  \big)
\end{align*}
are confidence intervals for the unknown parameter $F_0$ at the same levels given in the theorem. 
To prove this last assertion, we prove the following conditions that result sufficient: 
\begin{gather*}
P\Big( \hN{p_0}^{-1} (\lambda_0\bar{\mathcal{Y}} - h_d)  \geq F_0 \Big) 
\leq 1 - {\alpha}_+ \,;
\\
P\Big( \hN{p_0}^{-1} (\lambda_0\bar{\mathcal{Y}} + h_u) \leq F_0 \Big) 
\leq 1 -  {\alpha}_- \,.
\end{gather*}
Observe that, since the function $\hN{p_0}$ is invertible with continuous inverse (see \cite[Section~\sref{app:Spec_func}]{Ale20SUPP}), 
we get
\begin{gather*}
P\Big( \hN{p_0}^{-1} (\lambda_0\bar{\mathcal{Y}} - h_d)  \geq F_0 \Big) 
= P\Big( \bar{\mathcal{Y}} \geq \frac{\hN{p_0}(F_0) +h_d}{\lambda_0} \Big) 
\,;
\\
P\Big( \hN{p_0}^{-1} (\lambda_0\bar{\mathcal{Y}} + h_u) \leq F_0 \Big) 
= P\Big( \bar{\mathcal{Y}} \leq \frac{\hN{p_0}(F_0) -h_u}{\lambda_0} \Big) 
\,;
\end{gather*}
and hence the final result is a consequence of the 
following steps, that are proved in \cite[Section~\sref{sec:SMproof}]{Ale20SUPP}.
\begin{description}
\item[First step] the following two inequalities
\begin{gather*}
P\Big( \bar{\mathcal{Y}} \geq E(\bar{\mathcal{Y}})+\frac{h_d}{\lambda_0} \Big)  \leq 1 -  {\alpha}_+ ;
\\
P\Big( \bar{\mathcal{Y}} \leq E(\bar{\mathcal{Y}})-\frac{h_u}{\lambda_0} \Big)  \leq 1 -  {\alpha}_- 
\end{gather*}
are consequence of Chernoff bound inequalities;
\item[Second step]  the special function $\hN{p_0}$ is such that
\[
E(\bar{\mathcal{Y}}) = \frac{\hN{p_0}(F_0)}{\lambda_0}. \qedhere
\]
\end{description}
\end{proof}

\subsection{Connection with extreme value theory}\label{sec:extr_value}
The main result of this paper is based on the the fact that 
the random variables $(\bar{Y}_{r\,c} )_{r,c}$ are
independent, conditioned on $\boldsymbol{m}_c $, see Lemma~\ref{cor:Y_rc}. 
As discussed in Section~\ref{sec:stat_descr} and used in \cite[\seqref{eq:barYrcDef}]{Ale20SUPP}, these variables are given 
as the maximum of a random number of independent exponentially distributed random 
variables
\[
\bar{Y}_{r\,c} = 
\max_{\substack{o_1, \ldots,o_{m_{r\,c}}\colon R(c,o_j)=r\\ o_j \text{ different objects}}} \big( \bar{Y}(o_j,c) \big) ,
\]
A natural question is the relation of such considerations with the extreme value theory.
The well-known Fisher–Tippett–Gnedenko theorem \cite{Gned43} provides an asymptotic result, 
and it shows that, when $F_0\to\infty$, if there are sequences $a_{F_0}$ and $b_{F_0}$ such that
\((\bar{Y}_{r\,c}-a_{F_0})/b_{F_0}\) converges in law to a random variables $Z$, then
$Z$ must be Gumbel, Fr\'echet or Weibull (Type 1,2 or 3). 
In the proof of Theorem~\ref{thm:ris1}, we can recognize that
\[
E(e^{s(\bar{Y}_{r\,c}- E(\bar{Y}_{r\,c}))} ) = \prod_{j=1}^{m_{r\,c}} 
\frac{
e^{-\tfrac{s}{j\lambda_0}} } {
1 - \tfrac{s}{j\lambda_0}
}
\mathop{\longrightarrow}_{F_0\to\infty}
\Big(\Gamma(1-\tfrac{s}{\lambda_0}) e^{-\gamma \tfrac{s}{\lambda_0}} \Big) = E(e^{Z}),
\]
from which we can recognize that $Z$ has a Gumbell law. Since the Chernoff bounds on the mean
of such variables gives the same concentration inequalities as in Theorem~\ref{thm:ris1},
our result gives also the confidence interval based on the Chernoff bounds of the asymptotic distribution
based on the extreme value theory. 
In addition, note that \(E(e^{s(\bar{Y}_{r\,c}- E(\bar{Y}_{r\,c}))} ) \nearrow E(e^{Z})\), meaning
that the limit bounds is a analytic \emph{upper bound} for the 
concentration inequality,
that is the key point in the proof of Theorem~\ref{thm:ris1}.

\section{Analytical asymptotic discussion}\label{sec:analyticalDiscussion}
In this section we discuss the accuracy of the analytical approximation given in the main result to show
the appropriateness in this context.

\begin{figure}[tb]
\begin{center}
\centering
\fbox{\includegraphics[width= 0.45\textwidth]{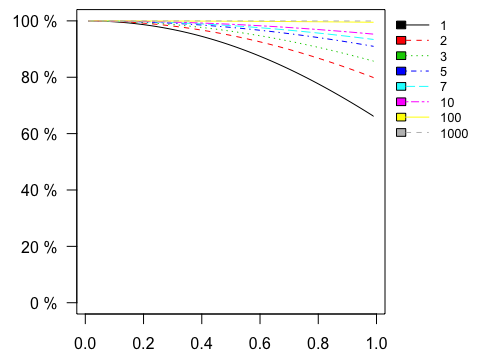}}
\fbox{\includegraphics[width= 0.45\textwidth]{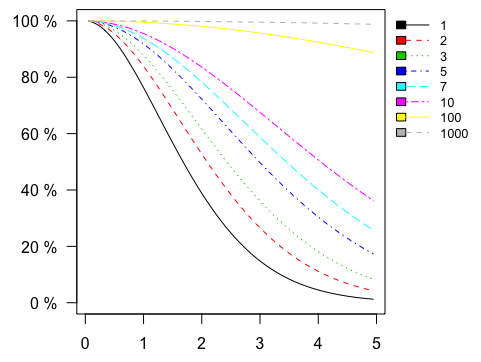}}
\end{center}
\caption{Ratio between the finite products and the series quantities given in
\eqref{eq:apprGamma}, for different values of $m_{r\,c}$ and $t$, expressed as
percentage of $\Gamma(1\mp t) e^{\mp \gamma t}$ given by
$\prod_1^{m_a} \frac{e^{\mp \frac{t}{j}}}{1 \mp \frac{t}{j}} $.
The different lines refer to different values of $m_{r\,c}$, given in the legend.
Left: percentage of approximation for $\Gamma(1-t) e^{-\gamma t}$,
$t\in (0,1)$. Right: percentage of approximation for $\Gamma(1+ t) e^{+\gamma t}$,
$t\in (0,5)$.} \label{fig:approx}
\end{figure}

The confidence intervals in this paper are based on the
uniform bounds given in the proofs of  Theorem~\ref{thm:ris1} with the following inequalities:
\begin{equation}\label{eq:apprGamma}
\begin{aligned}
&\text{for }{\alpha}_+: 
&&
\prod_{j=1}^{m_{r\, c}} \frac{e^{-\frac{t}{j}}}{1 - \frac{t}{j}} 
\leq
\Big(\prod_{j=1}^\infty \frac{e^{-\frac{t}{j}}}{1 - \frac{t}{j}} \Big) = 
\Gamma(1-t) e^{-\gamma t} 
, && t \in (0,1);
\\
&\text{for }{\alpha}_-: 
&& 
\prod_{j=1}^{m_{r\, c}} \frac{e^{\frac{t}{j}}}{1 + \frac{t}{j}} 
\leq
\Big(\prod_{j=1}^\infty \frac{e^{\frac{t}{j}}}{1 + \frac{t}{j}} \Big) = 
\Gamma(1+t) e^{\gamma t} 
, && t >0
.
\end{aligned}
\end{equation}
We recall that $m_{r\, c}$
is the (random) number of object assigned to register $r$ by the hash function $c$. 
In Figure~\ref{fig:approx} we underline that this approximation is good for small values of $t$ and
big $m_{r\, c}$. 
To show that the uniform bound in this paper does not affect significantly the Chernoff bounds, 
we compare for different values of $h_u$ and $h_d$:
\begin{equation}\label{eq:apprGamma2}
\begin{aligned}
&\text{for }{\alpha}_+: 
&&
\min_{t\in (0,1)} \Big( \prod_{c=1}^{{c_0}}\prod_{r=1}^{2^{{r}_0}} 
e^{-t h_d}
\prod_{j=1}^{m_{r c}} \frac{e^{-\frac{t}{j}}}{1 - \frac{t}{j}} 
\Big)
&&\text{vs.}&&
\Big( \Gamma(1-t_+) e^{-(\gamma +h_d) t_+} \Big)^{c_02^{{r}_0}} 
;
\\
&\text{for }{\alpha}_-: 
&& 
\min_{t>0} \Big( \prod_{c=1}^{{c_0}}\prod_{r=1}^{2^{{r}_0}} 
e^{-t h_u}
\prod_{j=1}^{m_{r c}} \frac{e^{\frac{t}{j}}}{1 + \frac{t}{j}} 
\Big)
&&\text{vs.}&&
\Big( \Gamma(1+t_-) e^{(\gamma -h_u) t_-} \Big)^{c_02^{{r}_0}} 
.
\end{aligned}
\end{equation}
\begin{figure}[tbhp]
\begin{center}
\includegraphics[width=0.95\linewidth]{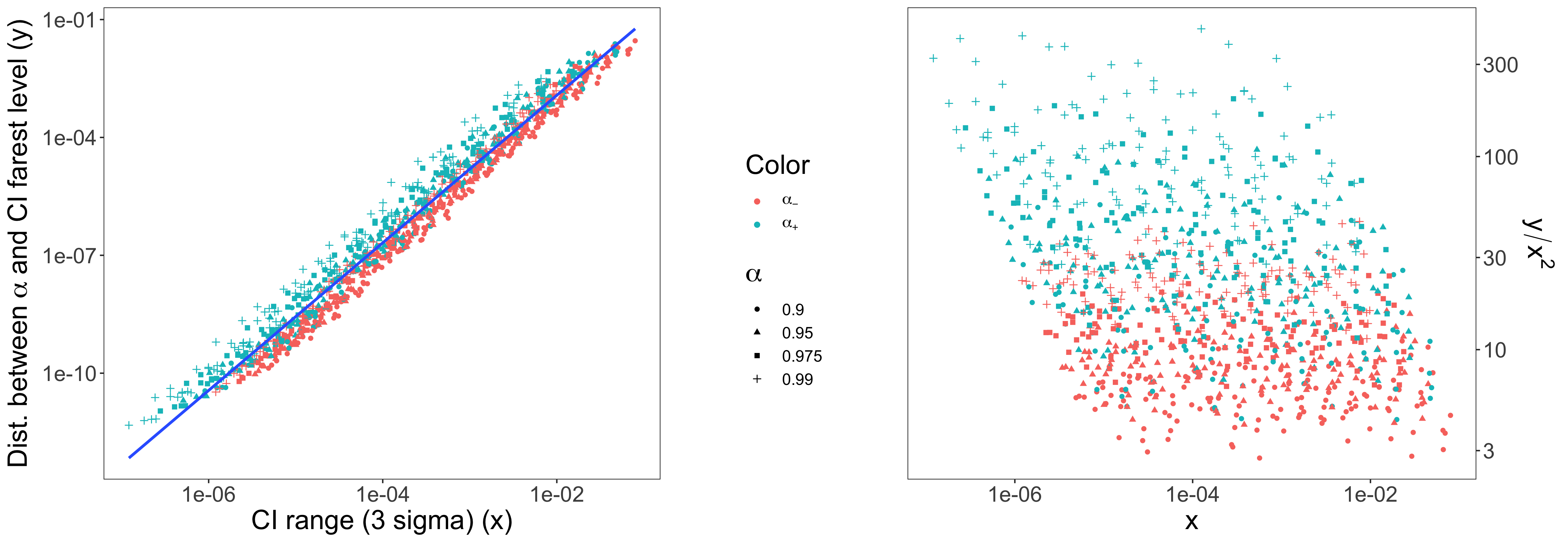}\\
\end{center}
\caption{Accuracy in the use of the analytical limit in \eqref{eq:apprGamma2} (MonteCarlo simulation of 
$\{\boldsymbol{m}_c , c= 1,\ldots,c_0\}$). Each point refers
to a different choice of ${\alpha}_+$ (light blue) or ${\alpha}_-$ (light red), $r_0 \in \{0,\ldots,4\}$, $c_0\in \{1,\ldots,4\}$
and $F_0 \in \{50,100,500,1000,5000,10000,50000,100000\}$.  
Left: linear dependence
in log-log scale ($y = 1.91 + 1.88 x $) between the precision in using the exact formula ($x$ is the length
of the $3\sigma$ confidence interval of ${A}_\pm$) and the accuracy of the estimation of ${\alpha}$ with gamma
function instead of the exact formula ($y$ is the distance between ${\alpha}$ calculated with the gamma function and
the farthest endpoint of the $3\sigma$ exact confidence interval). Rigth: dependence in log-log scale of $y/x^2$ 
with respect to $x$ as function of different ${\alpha}_{\pm}$.}\label{fig:cfrFinal}
\end{figure}
For $r_0 \in \{0,\ldots,4\}$, $c_0\in \{1,\ldots,4\}$, 
and $\alpha \in \{.9,.95,.975,.99\}$,
we choose the values of $h_u$ and $h_d$ for which 
\[
\Big( \Gamma(1-t_+) e^{-(\gamma +h_d) t_+} \Big)^{c_02^{{r}_0}} 
= 1-\alpha_\pm = 
\Big( \Gamma(1+t_-) e^{(\gamma -h_u) t_-} \Big)^{c_02^{{r}_0}} .
\]
Then, for any
$F_0 \in \{50,100,500,1000,5000,10000,50000,100000\}$, 
with a MonteCarlo procedure, we estimate the mean value and the standard deviation of the
random quantities
\[
{A}_- = \min_{t\in (0,1)} \Big( \prod_{c=1}^{{c_0}}\prod_{r=1}^{2^{{r}_0}} 
e^{-t x_d}
\prod_{j=1}^{m_{r c}} \frac{e^{-\frac{t}{j}}}{1 - \frac{t}{j}} 
\Big)
\qquad 
\text{ and }
\qquad 
{A}_+ = \min_{t>0} \Big( \prod_{c=1}^{{c_0}}\prod_{r=1}^{2^{{r}_0}} 
e^{-t x_u}
\prod_{j=1}^{m_{r c}} \frac{e^{\frac{t}{j}}}{1 + \frac{t}{j}} 
\Big)
\]
by simulating different values of the multinomial vectors $\{\boldsymbol{m}_c , c= 1,\ldots,c_0\}$. 
As expected, all the simulated quantities above result smaller than $1-\alpha$. Then, 
for each $r_0,c_0,\alpha,F_0$ we have built a $3\sigma$ confidence interval $[{\mathfrak{a}}_-^l,{\mathfrak{a}}_-^u]$ and $[{\mathfrak{a}}_+^l,{\mathfrak{a}}_+^u]$ 
for ${A}_-$ and ${A}_+$, respectively.
All the results are presented in Figure~\ref{fig:cfrFinal}. On the left-hand side, it is drawn 
the scatter-plot of 
\begin{align*}
&x = \text{range of confidence interval} &&= {\mathfrak{a}}_+^u -{\mathfrak{a}}_+^l &&\text{ (${\mathfrak{a}}_-^u -{\mathfrak{a}}_-^l$, respectively)};
\\
&y = \text{maximum imprecision} &&= {{\mathfrak{a}}}_+ - {\mathfrak{a}}_+^l &&\text{ (${{\mathfrak{a}}}_- -{\mathfrak{a}}_-^l$, respectively)};
\end{align*} 
which shows a good linear dependence in a log-log scale. As the linear coefficient is close to $2$,
on the right-hand side, the scatterplot of $y/x^2$ vs.\ $x$ confirms this scale of dependence, and
it suggests that the variability of the constant depends mainly on ${\alpha}$, 
firstly on the choice of the sign ($\alpha_+$ or $\alpha_-$), and then on its value.

A finer analysis shows that, when $F_0 \geq 500$, the maximum imprecision 
is less than $0.00683$ (with $r_0=4$, $c=1$, $p_-= 0.1$, $N_0 = 500$), 
becoming less than $6.7 \cdot 10^{-5}$ for $F_0 \geq 50000$ (again, $r_0=4$, $c=1$, $p_-= 0.1$ but $N_0 = 50000$). 
In other words, the uniform bounds given in \eqref{eq:apprGamma2} appear adequate in this context.

\section{Application on a real data-stream}\label{sec:realDataAnalysis}
We test the algorithms described above on Twitter data (with unique user IDs $F_0 = 454,176$)
and on an anonymized real time data stream, made by 
$196,432,300$ objects, of which $F_0 = 1,407,593$ distinct. 

\begin{figure}[th]
\begin{center}
\includegraphics[width=0.45\linewidth]{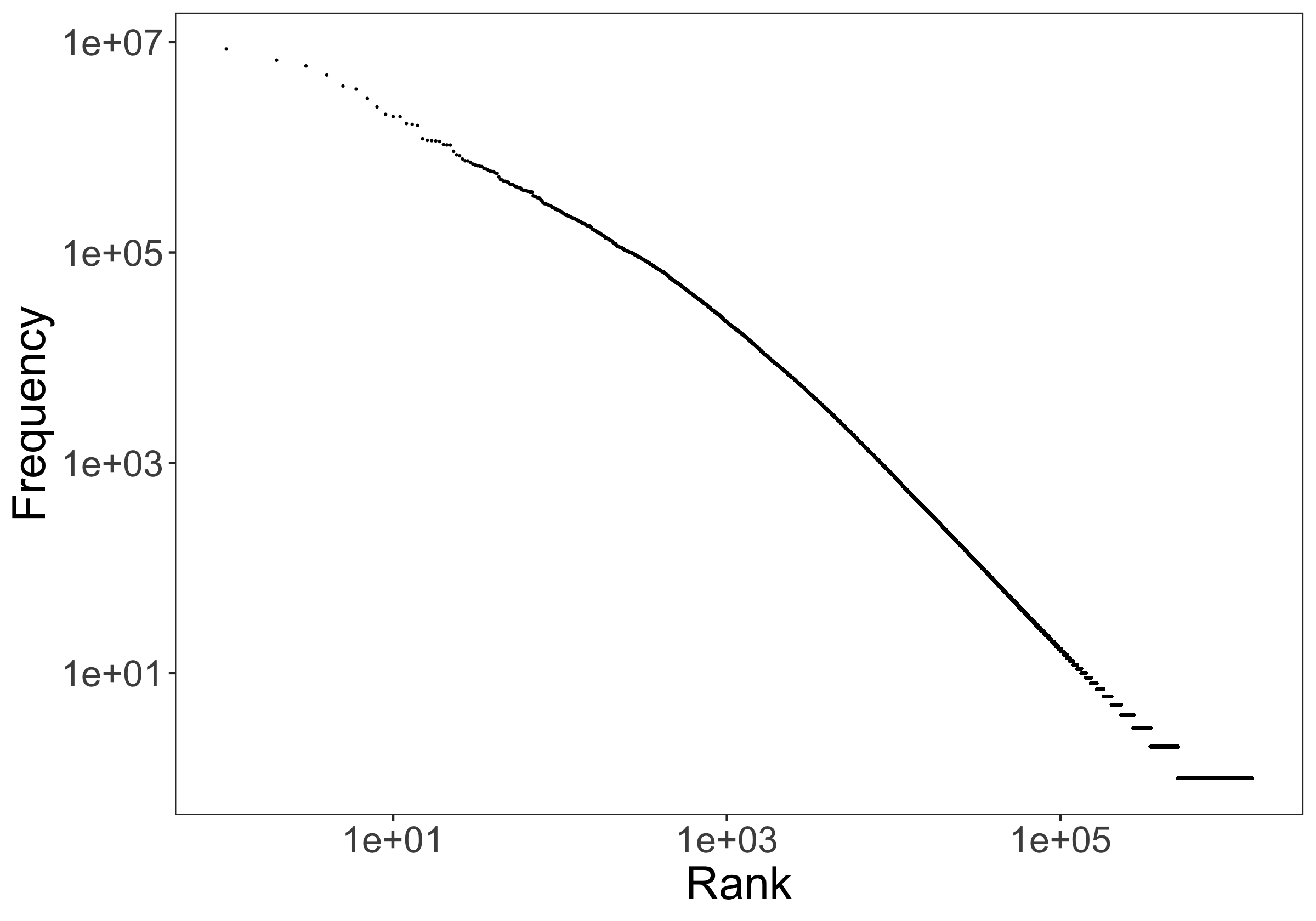}\\
\end{center}
\caption{Frequency rank plot of the frequency count of the $F_0 = 1,407,593$ distinct objects in the real data stream. 
The log-log linear plot indicates the good fit to the power law distribution.}\label{fig:powerLaw}
\end{figure}

The distribution of the occurrences of the second bigger database
may be seen as a power law distribution, as shown by
the log-log frequency rank plot (see Figure~\ref{fig:powerLaw}).

The data are divided into compressed files ($100$ for Twitter data and $1,000$ for real time), 
and analyzed with Apache Spark on R.
The $\text{SHA}_{256}$ function $\texttt{sha2(Id, 256)}$ has been applied to each object, and
the $256$-bits output has been divided into $4$ equal parts, each of one being
certified to be a sequence of i.i.d.\ Bernoulli random variables 
(see \cite{HTML1,HTML3,HTML2}). With such a division, we analyze our data-stream with
$c_0 = 4$ hash functions. 
Moreover, since Spark codes $\texttt{sha2}$ output as a hexadecimal string, we used the first
character ($4$ bits) to define $r_0 = 4$, so that we have $a_0 = 4 \cdot 2^4 = 64$ registers where
we store the values of $\mathbb{Y}$ and $\mathbb{Z}$ during the streaming algorithm, and the last
$2$ characters to define $z_0 = 8$, noticing that the remaining $13$ characters
($52$ bits) are sufficient for the definition of $\mathbb{X}$ in this application.
These $13$ hexadecimal characters are converted into a binary string
and the number of leading $0$s are computed with $\texttt{52-length(binary string)}$.

\begin{figure}[tbhp]
\begin{center}
\centering
\fbox{\includegraphics[width= 0.45\textwidth]{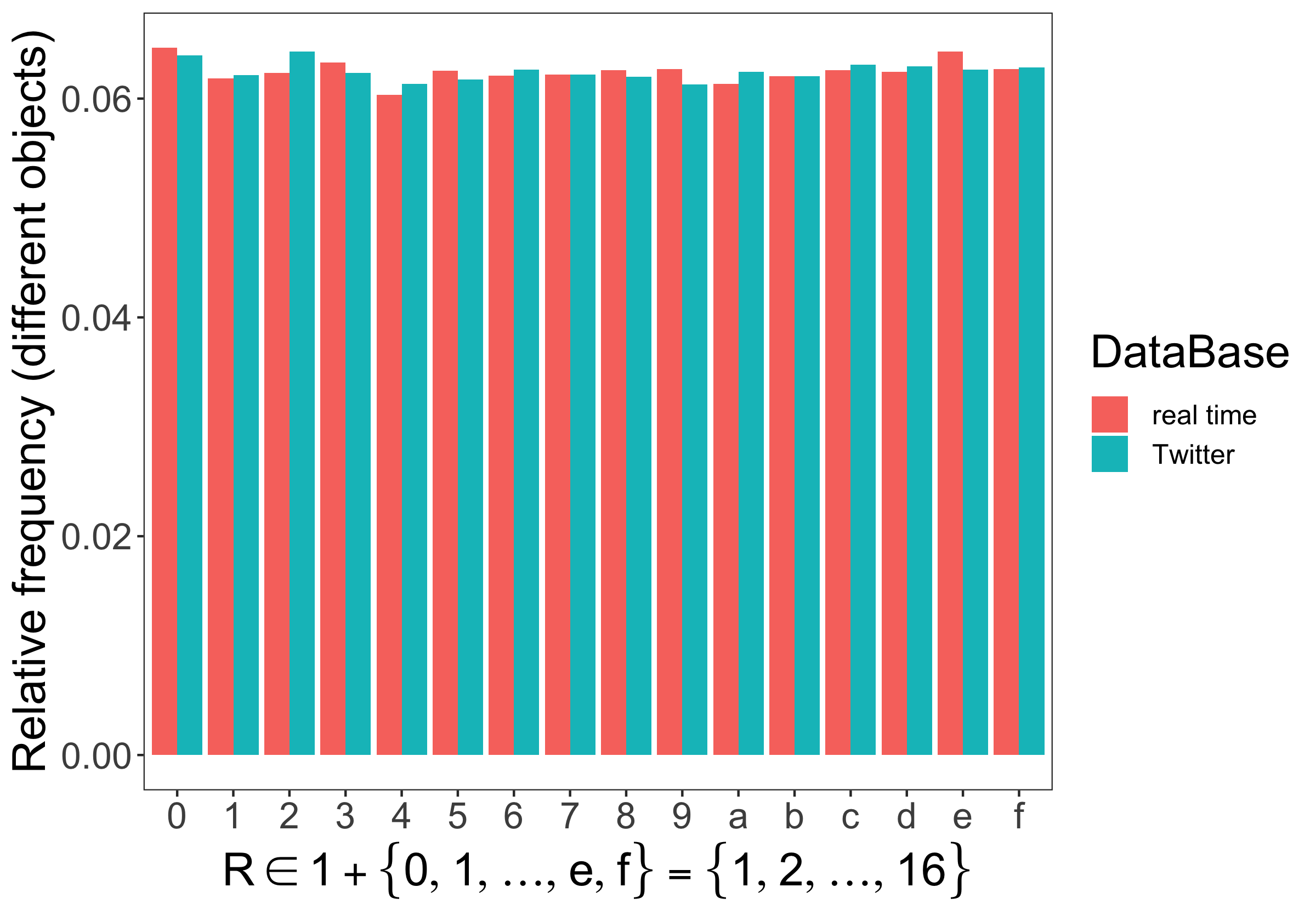}}
\fbox{\includegraphics[width= 0.45\textwidth]{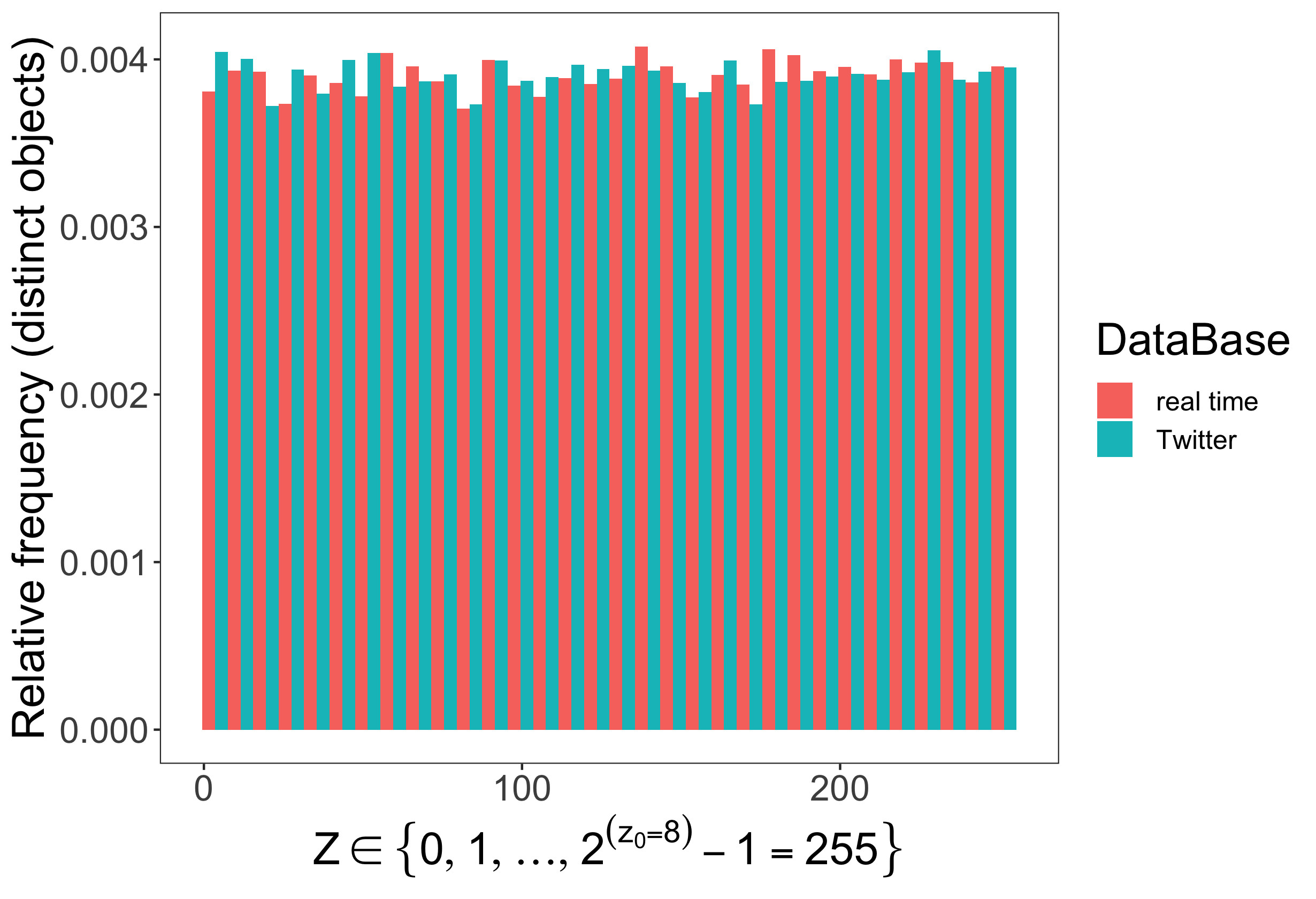}}
\fbox{\includegraphics[width= 0.45\textwidth]{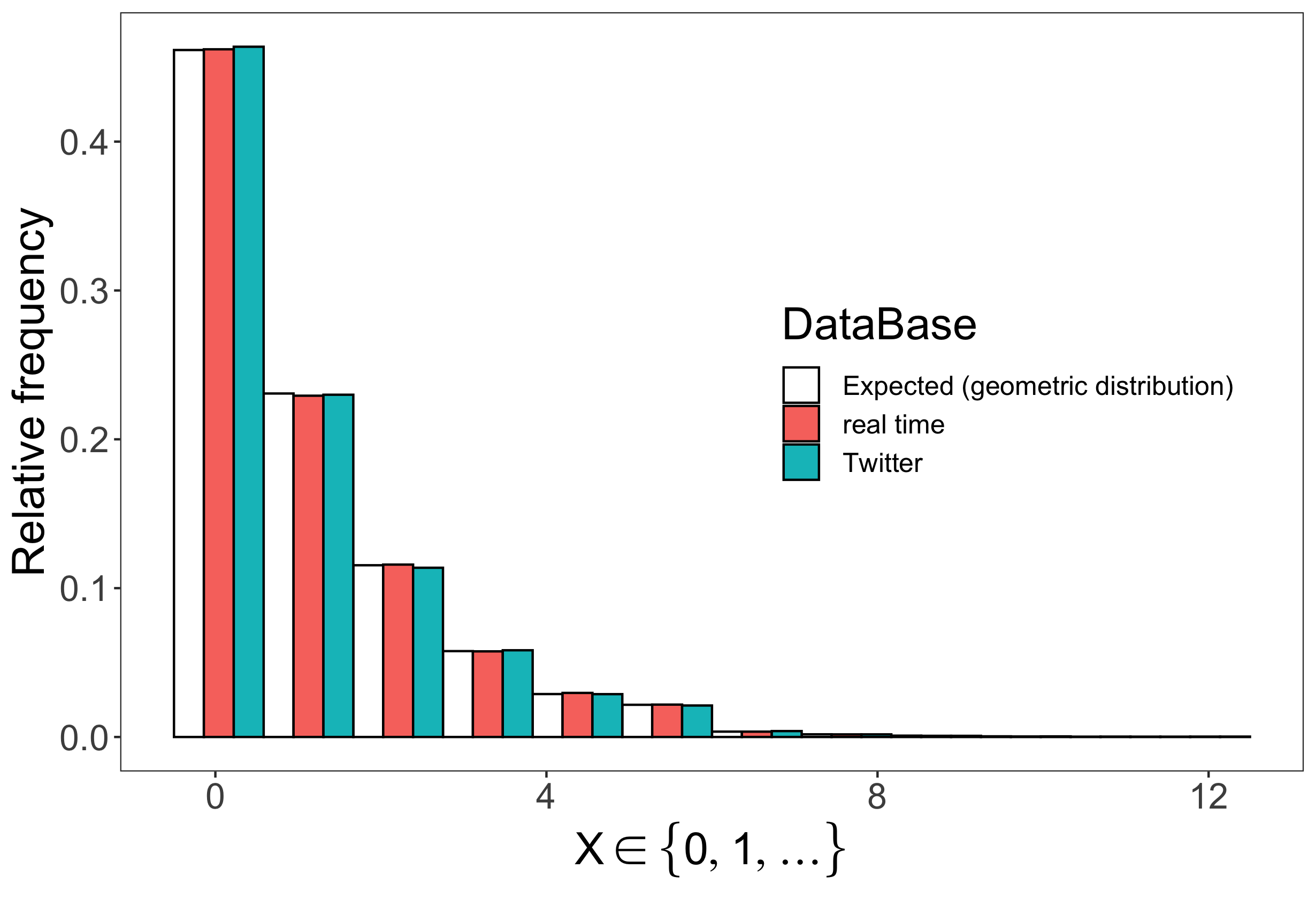}}
\fbox{\includegraphics[width= 0.45\textwidth]{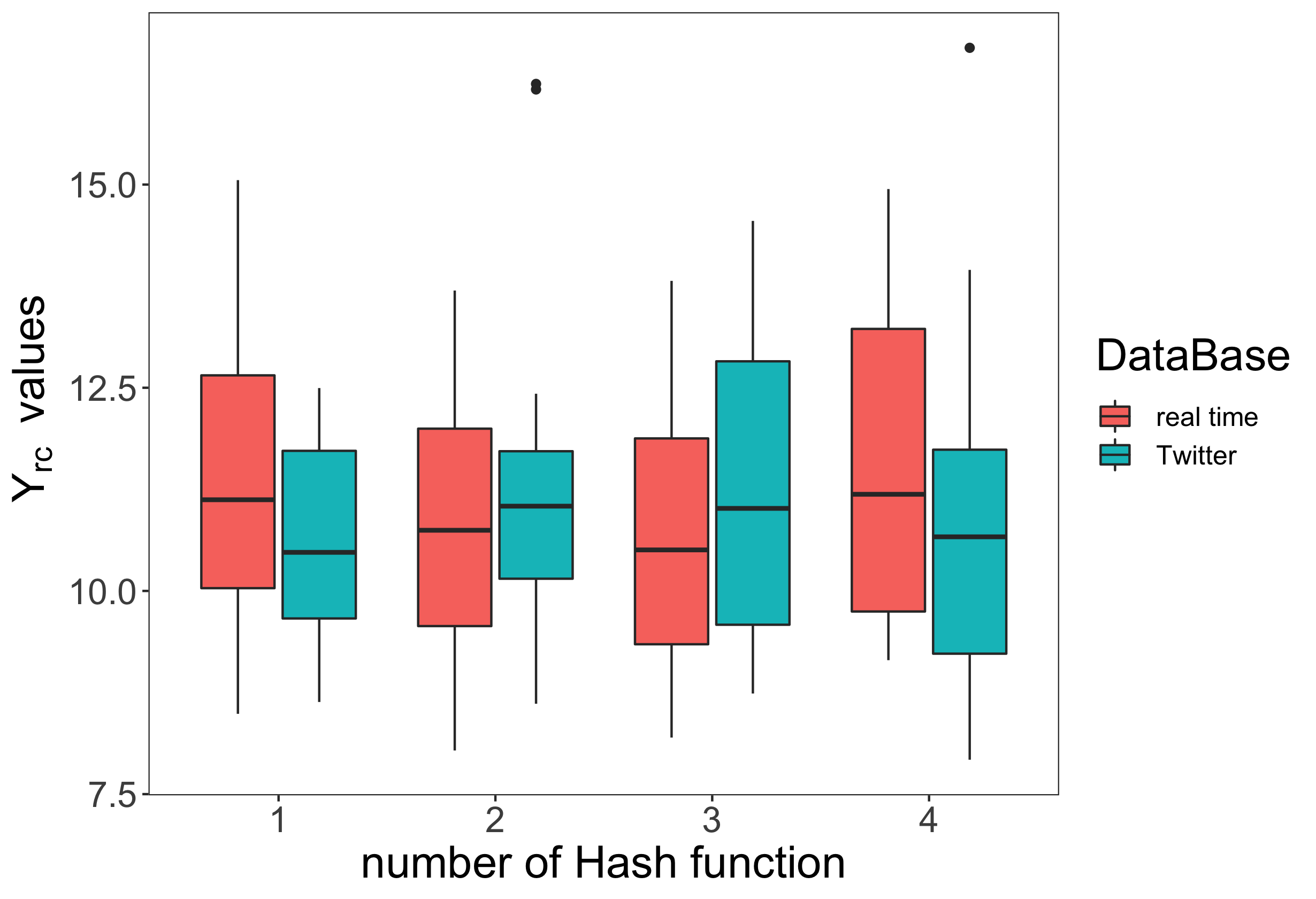}}
\end{center}
\caption{Randomness of hash generations in the first file (out of $1,000$). 
Top left: check of theoretical uniform distribution of the quantity $R$ evaluated only on the different objects
(domain made by $2^4=16$ possible different outcomes).
Top right: check of theoretical uniform uniform distribution of the quantity $Z$ evaluated only on the different objects
(domain made by $2^8=256$ possible different outcomes).
Bottom left: check of theoretical distribution of the quantity $X$ evaluated only on the different objects 
compared with the expected geometric distribution.
Bottom right: spread of the quantities $\mathbb{Y}$ divided by the hash function: 
each boxplot groups $\{{Y}_{r\,c}, r=1,\ldots,2^4\}$ for different $c=1,2,3,4.$} \label{fig:uniform}
\end{figure}

\subsection*{Goodness of fit of statistical distributions}
Before giving the overall results, we analyze the results of a single file for Twitter and ``real time'' datasets.
The stream data $\{o_1,o_2, \ldots\}$ is made by $49,999$ objects (resp.~$196,433$), 
made by $15,999$ different repeated objects (resp.~$18,094$).
Each object is signed with $4$ hash functions. We check the uniform distribution on the distinct objects
for the random values of $r \in \{1,\ldots,2^{r_0}= 16\}$ 
(Twitter: $\chi^2 = 17.608$, $df = 15$, $\text{p-value} = 0.2838$)
(real time: $\chi^2 = 10.597$, $df = 15$, $\text{p-value} = 0.7808$)
and of $z \in \{1, \ldots, 2^{z_0}= 256\}$ (Twitter: $\chi^2 = 298.91$, $df = 255$, $\text{p-value} = 0.03062$, 
read data: $\chi^2 = 286.51$, $df = 255$, $\text{p-value} = 0.08521$),
and of the geometric distribution of $X$ (Twitter: $\chi^2 = 8.7522$, $df = 12$, $\text{p-value} = 0.7239$, 
real time: $\chi^2 = 7.1689$, $df = 12$, $\text{p-value} = 0.8463$).
We plot the corresponding histograms in Figure~\ref{fig:uniform}, together with the boxplots of the
registers $\mathbb{Y} = Y_{r\,c}$ grouped by $c$, the hash key 
(ANOVA test: Twitter $F_{3,60} = 0.437$, $\text{p-value} = 0.728$, real time $F_{3,60} = 1.095$, $\text{p-value} = 0.358$).

\subsection*{Accuracy of the algorithm}
\begin{figure}[tbhp]
\begin{center}
\centering
\fbox{\includegraphics[width= 0.3\textwidth]{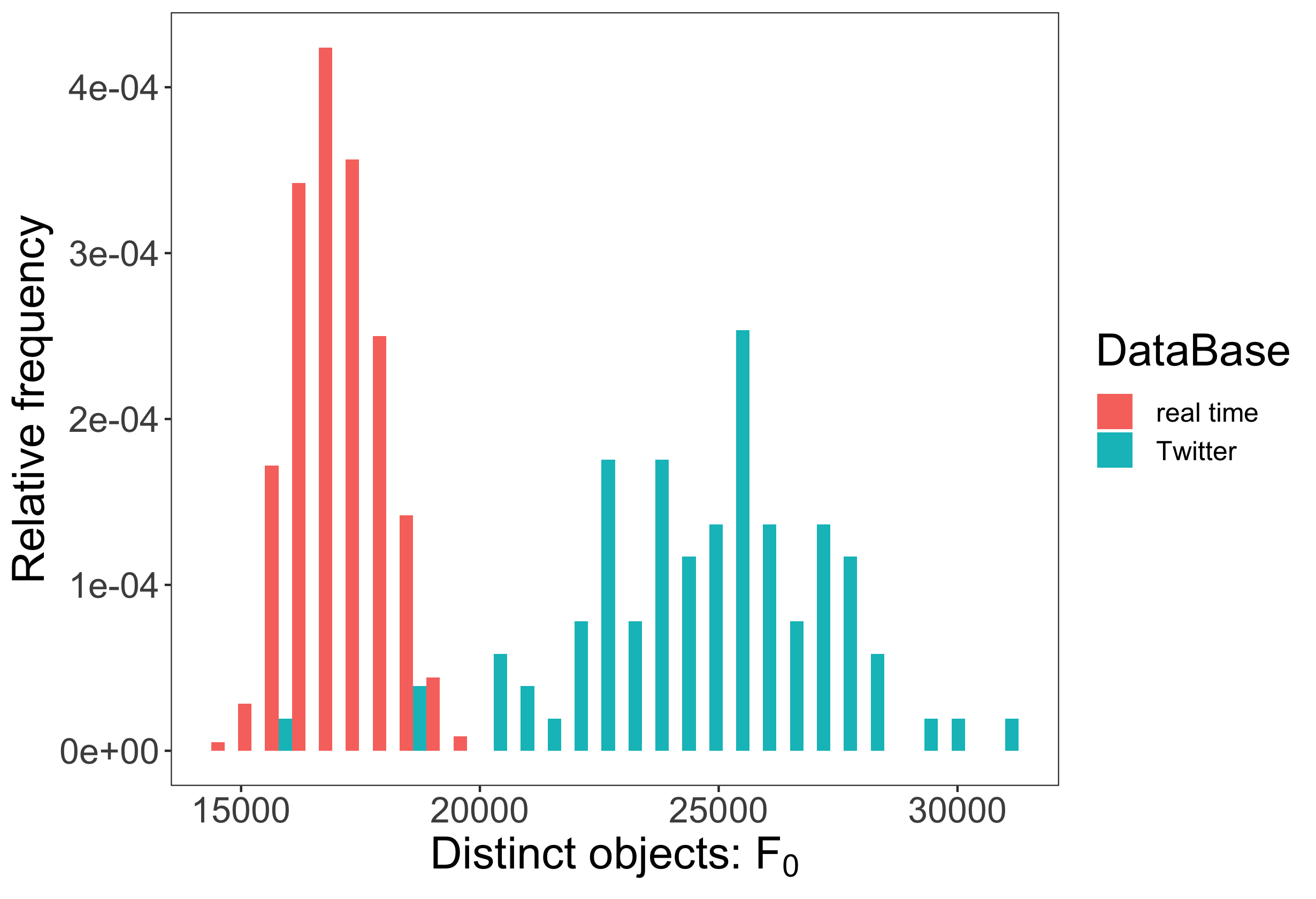}}
\fbox{\includegraphics[width= 0.3\textwidth]{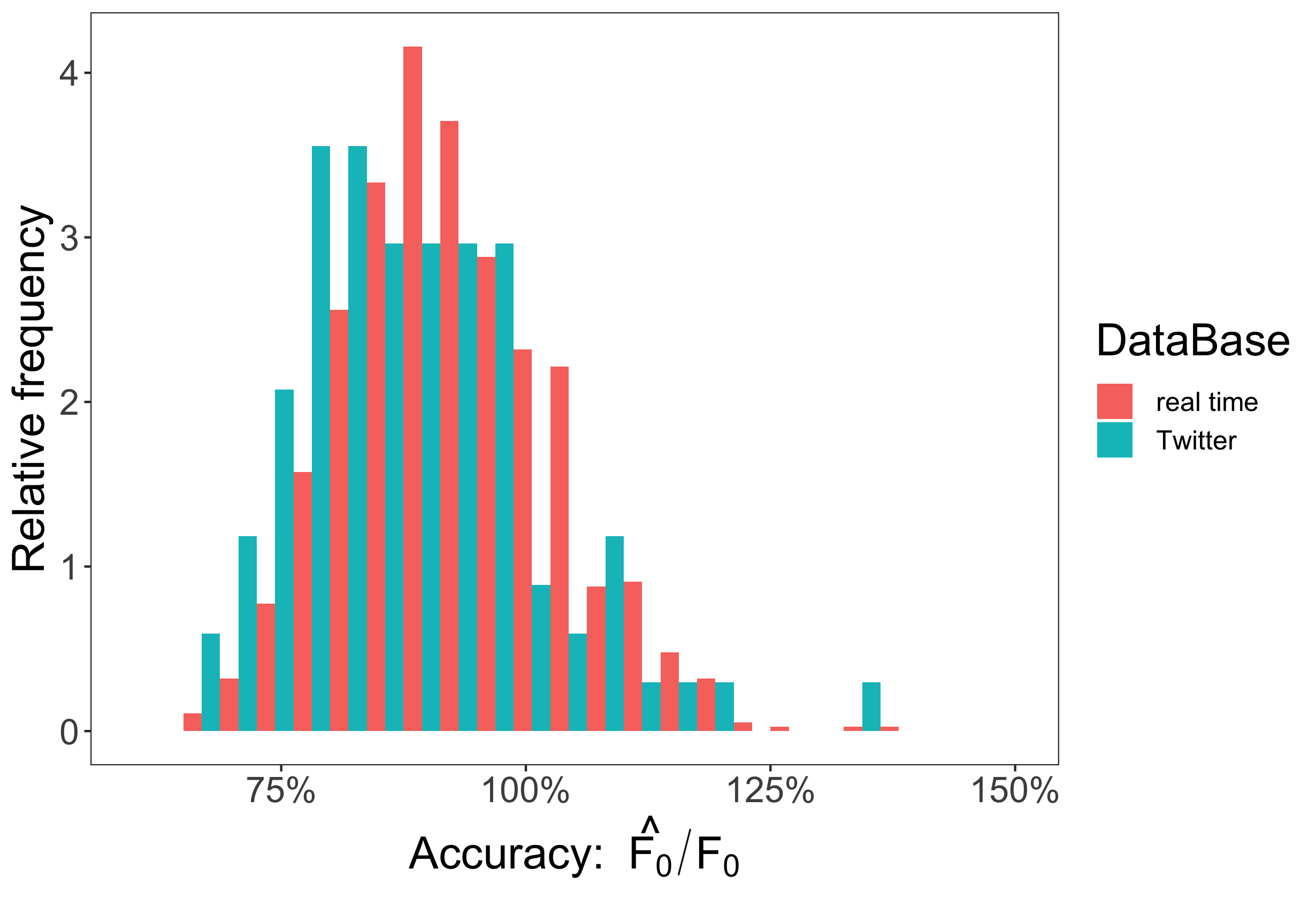}}
\fbox{\includegraphics[width= 0.3\textwidth]{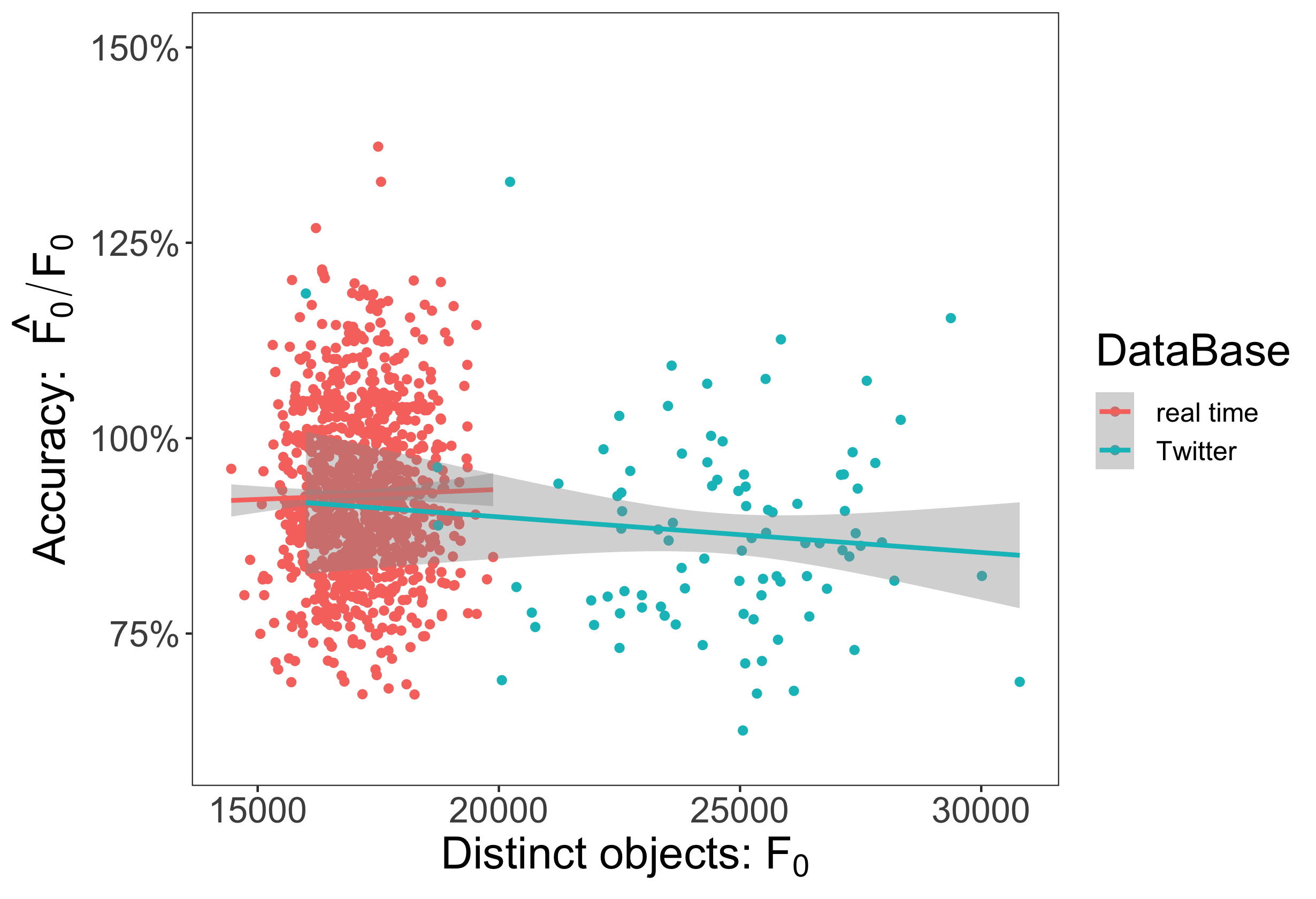}}
\fbox{\includegraphics[width= 0.45\textwidth]{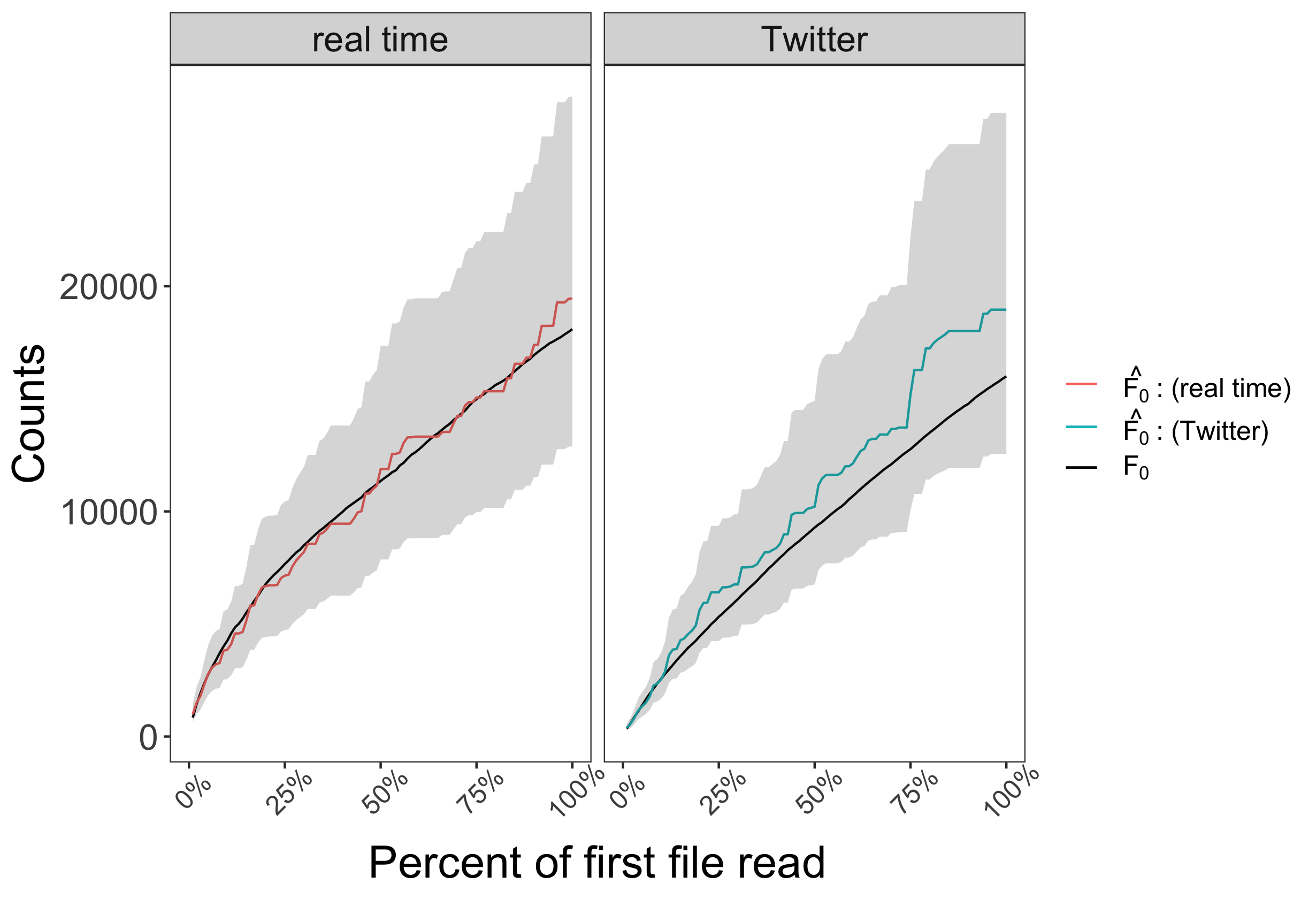}}
\fbox{\includegraphics[width= 0.45\textwidth]{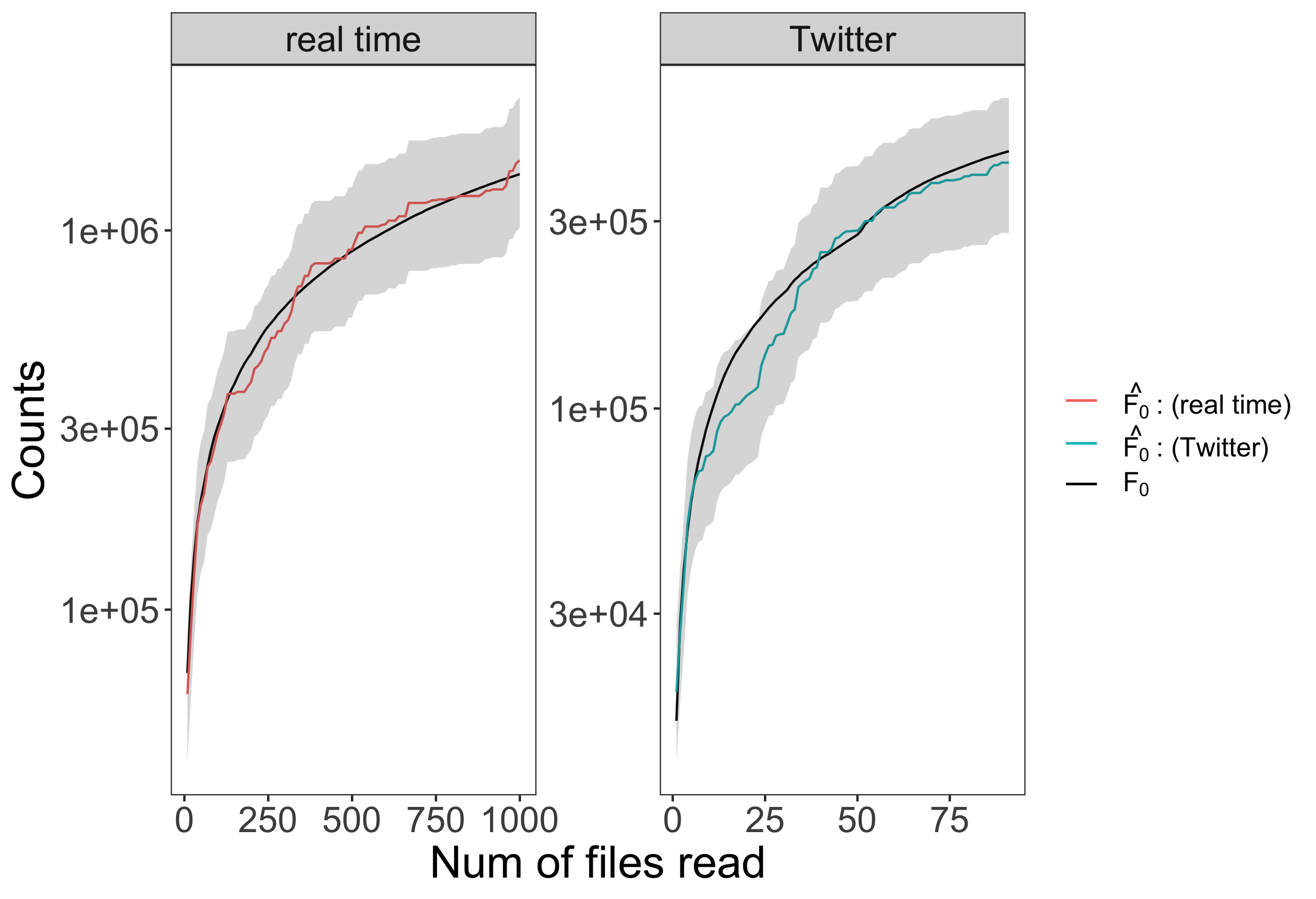}}
\end{center}
\caption{Analysis of accuracy of the estimations.
Top-left: histogram of the number $F_0$ of distinct object in each file.
Top-center: histogram of the percent accuracy $\hat{F_0}/F_0$ of 
the estimates of $F_0$ made on each file. 
Top-right: scatterplot of the relative accuracy $\hat{F_0}/F_0$ vs.\ the number $F_0$ of distinct object in each file.
Bottom-left: evolution of the confidence interval for $F_0$ during the analysis of the data in the first file (in black: true value of $F_0$).
Bottom-right: evolution of the confidence interval for $F_0$ during the analysis of all the streaming data (in black: true value of $F_0$).} 
\label{fig:countAlgo}
\end{figure}

We then analyze each of the compressed files, that contains a different value of distinct object $F_0$.
The distribution of the true $F_0$ is plotted in Figure~\ref{fig:countAlgo} (top-left).
For each of this file, we also estimate $F_0$ with $\hat{F_0} = \hN{p_0}^{-1}  (\lambda_0\mathcal{Y} ) $,
and we compute the relative accuracy of each estimation with $\hat{F_0} /F_0$.
The distribution of the relative accuracy is plotted in Figure~\ref{fig:countAlgo} (top-center) for both the databases.
In Figure~\ref{fig:countAlgo} (top-rigth), the scatterplot of the accuracy $\hat{F_0} /F_0$ vs.\ $F_0$ 
shows that there is not association between these two variables
(Twitter $R^2 = 0.008839$, $\text{p-value} = 0.3754$, real time $R^2 = 0.0004698$, $\text{p-value} = 0.494$).

Finally, we analyze the data sequentially as a data stream. 
We check that the $90\%$ confidence interval is consistent all along the process. In
Figure~\ref{fig:countAlgo} it is shown the evolution of the confidence interval at the beginning of the stream
(during the first file, bottom left) and its consistency as the number of files increases ($\log$-scale, bottom left). As expected, 
the cold-start effect is mitigated since the approximation is not made on asymptotic properties.

\section{Theoretical resolution of computational aspects}\label{sec:CompAsps}
We recall that we build confidence intervals for $F_0$ based on the output $\mathbb{Y}$
of Algorithm~\ref{alg:querying_data}.
For example, 
Algorithm~\ref{alg:conf_int} shows how to compute the confidence interval of the form $(0, \mathrm{upper})$
and it faces two nonlinear problems.
Analogous procedures can be used to compute confidence intervals of other forms.
The key computational point is the
necessity of numerically solving some nonlinear equations that involve mathematical special functions.

In the following sections, we state the relevant inequalities that can be used to find the root of $f(x)=0$ in our context,
with the Halley's method \cite{Scavo95}.
This iterative method is given by
\[
x_{n+1} = x_n - \frac {2 f(x_n) f'(x_n)} {2 {\big(f'(x_n)\big)}^2 - f(x_n) f''(x_n)},
\]
it is essentially the Newton method applied to the function $g(x) = \tfrac{f(x)} {\sqrt{|f'(x)|}}$,
and it achieves a cubic rate of convergence in the neighborhood of the solution, see \cite{Alef81}.

In addition, we give accurate lower and upper bounds for the solution, that can be shown to be contained in
the basin of attraction of the solution. Note that these bounds can be used also with a much simpler and
robust bisection method, which has as the counterpart a linear rate of convergence. 

\subsection{The problem $\psi(x) - y = 0$}
We recall here that the \emph{digamma function} $\psi:(0,\infty)\to \mathbb{R}$ is defined as
the logarithmic derivative of the $\Gamma$ function, see \cite[\S 6.3]{AbsSteg64}, and it satisfies the
relation
\begin{equation}\label{eq:psiDef1}
\psi(x+1) = \psi(x)+\tfrac{1}{x} .
\end{equation}
In addition $\psi$ is 
a strictly monotone, concave function, with $\lim_{t\to 0^+}\psi(t)=-\infty$,
$\psi(1) = -\gamma$ and $\psi(t)=\log(t) +o(1)$ when $t\to\infty$ (see, for example, \cite{DS2016}). 
Finally, it is implemented in all the recent math packages together
with its first and second derivative functions $\psi_1$ and $\psi_2$. 

As shown in Section~\sref{sec:boundPsi}, we have
\begin{equation}\label{eq:psiAppr}
\ln (x-\tfrac{1}{2}) < y < \ln (x) , \qquad e^{y} < x < e^{y} + \tfrac{1}{2}, 
\qquad \forall x > \tfrac{1}{2}, \forall y= \psi(x).
\end{equation}

\subsection{The problem $\hN{p}(x)-y = 0$}
First note that $\hN{p}(x), \hN{p}'(x)$ and $\hN{p}''(x)$ may be computed with with arbitrary precision,
because of  \seqref{eq:PHarm3} and \seqref{eq:PHarmDer} and the fact that a quad-double precision 
algorithm to calculate Lerch's transcendent of real arguments have been already developed, see
\cite{AKSENOV20031}.

For $p\in(0,1)$, as shown in \cite[Section~\sref{sec:boundHnP}]{Ale20SUPP}, we have
\begin{equation}\label{eq:center_eqAll}
\frac{{e^{y-\gamma}}}{p} - \frac{1}{2} 
\geq
x \geq
\begin{cases}
\frac{{e^{y-\gamma}}}{p} - e +\tfrac{1}{\ln(1-p)} & 
\text{if }{y} > \log \big( \gamma + p(\tfrac{1}{2} -\tfrac{1}{(e-1)\ln(1-p)}) \big); \\
{e^{y-\gamma}} - 1 & 
\text{otherwise}.
\end{cases}
\end{equation}

\subsection{The problem $y= (x-\gamma) t(x) - \ln \Gamma(1+t(x)) $, where $t(x)= \psi^{-1}(x-\gamma)-1$}
Note that, if $g(x)= (x-\gamma) t(x) - \ln \Gamma(1+t(x))$, then
\begin{equation}\label{eq:p-der}
g'(x) =  t(x) + t'(x) (x-\gamma - \psi(1+t(x)) ) = t(x),
\end{equation}
since, by definition of $t(x)$, $\psi(1+t(x)) = x-\gamma$. Then  
the formula of the derivative of the inverse function gives
\[
g''(x) =  t'(x) =  
\frac{1}{\psi_1(\psi^{-1}(x-\gamma))} = \frac{1}{\psi_1(1+t(x))} .
\]
As shown in \cite[Section~\sref{sec:boundAm}]{Ale20SUPP}, we have
\begin{equation}\label{eq:boundsAmAll}
\begin{aligned}
\sqrt{\frac{1}{50}  y}
&<
x 
< 
\pi\sqrt{\frac{2}{3}  y} , && \text{if $y <3$;}
\\
\frac{2}{3}  \Big( \log \Big( y +\frac{1}{2} \Big) + \gamma \Big) 
&< x < 2 \Big(\log \Big(\frac{4}{3}y+1\Big) + \gamma\Big) 
, && \text{if $y \geq 3$.}
\end{aligned}
\end{equation}

\subsection{The problem $y= (x+\gamma) t(x) - \ln \Gamma(1-t(x)) $, where $t(x)= 1-\psi^{-1}(-x-\gamma)$}
Note that, if $g(x)= (x+\gamma) t(x) - \ln \Gamma(1-t(x)) $, then
\begin{equation}\label{eq:p+der}
g'(x) =  t(x) + t'(x) (x+\gamma - \psi(1-t(x)) ) = t(x),
\end{equation}
since, by definition of $t(x)$, $\psi(1-t(x)) = -x-\gamma$. Then  
the formula of the derivative of the inverse function gives
\[
g''(x) =  t'(x) =  
\frac{1}{\psi_1(\psi^{-1}(-x-\gamma))} = \frac{1}{\psi_1(1-t(x))} .
\]

As shown in \cite[Section~\sref{sec:boundAp}]{Ale20SUPP}, we have
\begin{equation}\label{eq:boundsApAll}
\max\Big( - \ln (1 - C) - \gamma , \frac{\pi^2}{6} C \Big) 
<
x 
< 
2 \sqrt{(y+1)^2-1},
\end{equation}
where
\[
C = \sqrt{ 1 - \frac{-({\tfrac{y}{2}}-{\tfrac{6 + \pi^2}{12}}) + \sqrt{({\tfrac{y}{2}}-{\tfrac{6 + \pi^2}{12}})^2 +4{\tfrac{18 - \pi^2}{12}}} }{ 2 {\tfrac{18 - \pi^2}{12}}}} \in (0,1).
\]

\subsection{Minimum $\log$-length interval}
In this section, we show how to numerically compute the minimum length interval, in $\log$-scale, 
for a given confidence $\alpha$, based on the inequalities given in 
Theorem~\ref{thm:ris1}. 
The probem is set as follows: given $\alpha \in (0,1)$, $r_0 \geq 0$, $c_0 \geq 1$, 
we want to solve the nonlinear minimization problem:
\begin{gather*}
\min(h_d+h_u)
\intertext{subject to }
\begin{cases}
{\alpha}_+  = 1 - \exp\Big( -2^{{r}_0}c_0 \big[ (h_d+\gamma) t_+ - \ln \Gamma(1-t_+) \big]\Big), \qquad t_+ = 1-\psi^{-1}(-h_d-\gamma);
\\
{\alpha}_-  = 1 - \exp\Big( -2^{{r}_0}c_0 \big[  (h_u-\gamma)  t_- - \ln \Gamma(1+t_-) \big]\Big), \qquad t_- = \psi^{-1}(h_u-\gamma)-1;
\\
{\alpha}_+  + {\alpha}_- \geq 1 + \alpha;
\\
h_d,h_u\geq 0.
\end{cases}
\end{gather*}
The two values ${\alpha}_+$ and ${\alpha}_-$ are monotone functions of $h_d$ and $h_u$, respectively, as 
a consequence of \eqref{eq:p+der} and \eqref{eq:p-der}. As a consequence, the minimum is attained when ${\alpha}_+  + {\alpha}_- = 1 + \alpha$.
Then, 
if we set $x = 1 - {\alpha}_+$, 
the problem above may be rewritten
in terms of $x$: given ${\alpha} \in (0,1)$ and $a_0 = 2^{{r}_0}c_0 \in \{1,2,\ldots\}$, find
\begin{gather*}
\min (g(x)) =
\min \big(y_+^{-1}(-\tfrac{\log x}{a_0})+y_-^{-1}(-\tfrac{\log(1-{\alpha}-x)}{a_0}) \big) 
\intertext{subject to }
\begin{cases}
y_+(h)  = (h +\gamma) t_+ - \ln \Gamma(1- t_+), 
\qquad t_+ = 1-\psi^{-1}(-h-\gamma);
\\
y_-(h)  = (h -\gamma) t_- - \ln \Gamma(1+ t_-), 
\qquad t_- = \psi^{-1}(h-\gamma)-1;
\\
0 \leq x \leq 1-{\alpha}.
\end{cases}
\end{gather*}
%
%
Differentiating $g$ with respect to $x$, since $y_{\pm}'(h) = t_\pm(h)$ by \eqref{eq:p+der} and \eqref{eq:p-der}, we obtain,
\begin{equation*}
g'(x) = - \frac{1}{a_0x}\frac{1}{t_+ \Big(y_+^{-1} \Big(-\frac{\log x}{a_0} \Big)\Big)}
+ \frac{1}{a_0(1-{\alpha} - x)}\frac{1}{t_- \Big(y_-^{-1} \Big(-\frac{\log (1-{\alpha} - x)}{a_0} \Big)\Big)}
\end{equation*}
which is null when the following equation is zero
\begin{equation*}
f(x) = {x}{t_+ \Big(y_+^{-1} \Big(-\frac{\log x}{a_0} \Big)\Big)}
- {(1-{\alpha} - x)}{t_- \Big(y_-^{-1} \Big(-\frac{\log (1-{\alpha} - x)}{a_0} \Big)\Big)}
\end{equation*}
Call 
\begin{equation*}
\hat t_+ = \hat t_+ (x) = t_+ \Big(y_+^{-1} \Big(-\frac{\log x}{a_0} \Big)\Big), \quad
\hat t_- = \hat t_- (x) = t_- \Big(y_-^{-1} \Big(-\frac{\log (1-{\alpha} - x)}{a_0} \Big)\Big), 
\end{equation*}
$\psi_1(x) = d\frac{\psi(x)}{dx}$ and 
$\psi_2(x) = d\frac{\psi_1(x)}{dx}$,
then
\[
d\frac{{\hat{t}_+}(x)}{dx} = -\frac{1}{a_0x}\frac{1}{\hat t_+ \psi_1(1 - \hat t_+)}, \qquad
d\frac{{\hat{t}_-}(x)}{dx}  = +\frac{1}{a_0(1-{\alpha} - x)}\frac{1}{\hat t_- \psi_1(1 + \hat t_-)}.
\]
The problem is then to find the solution for the nonlinear problem $f(x)=0$ that may be solved 
with the Halley's method that involves the problems seen above, noticing that
\begin{align*}
f(x) & =  {x}{\hat{t}_+}
- {(1-{\alpha} - x)}{\hat{t}_- },
\\
f'(x) & = 
\hat{t}_+ - \frac{1}{a_0 \hat{t}_+ \psi_1(1 - \hat{t}_+)}
+ 
\hat{t}_- - \frac{1}{a_0 \hat{t}_- \psi_1(1 + \hat{t}_-)}
\\
f''(x) & = 
t'_+ \Big( 1 + \frac{\psi_1(1 - \hat{t}_+)-\hat{t}_+\psi_2(1 - \hat{t}_+)}{a_0(\hat{t}_+ \psi_1(1 - \hat{t}_+))^2}\Big)
\\
& \qquad +
t'_- \Big( 1 + \frac{\psi_1(1 + \hat{t}_-)+\hat{t}_+\psi_2(1 + \hat{t}_-)}{a_0(\hat{t}_- \psi_1(1 + \hat{t}_-))^2}\Big).
\end{align*}
and that a good starting point is given by $x_0 = \tfrac{1-{\alpha}}{2}$.

\section{Conclusions}
In this paper, we provide analytical confidence intervals for the number $F_0$ of distinct elements in data 
streams by analyzing a class of FMa. While the major concern of the state of the art is algorithm's complexity, 
here the new mathematical-statistical approach permits a extensive analysis of such classes of algorithms.
The HyperLogLog data structure (called $\mathbb{X}$ in this paper) is enriched with a 
new data matrix of fized size ($\mathbb{Z}$) that helps to bound uniformly the estimators during the
querying counting phase. In this phase, the Chernoff bounds may be applied analytically and gives
asymptotically efficient estimators that are related to the extreme value theory.
In addition, the relation \(
E(\bar{\mathcal{Y}}) = \frac{\hN{p_0}(F_0)}{\lambda_0}
\) introduces a new class of special functions $\hN{p_0}$ used to find the confidence interval.

Since the new theoretical results are based on some analytical, computational and numerical assumptions,
we have shown that these assumptions are always satisfied in real situations. First, the 
analytical asymptotic approximation made on Chernoff bounds is shown to be 
irrelevant when $F_0$ is large.
Then, statistical assumptions on the distributions of the quantities of interests are shown to be
satisfied on a real dataset and the accuracy of the methodology is provided. 
Finally, the computational solution of the problems related to 
the new special functions is solved by showing the basins of attraction for a Newton 
based method with cubic rate of convergence.

\newpage

\appendix 
\section*{\appendixname}
\numberwithin{equation}{section}

In this document we collect some technical results useful for \cite{Ale20}.
Therefore, the notation and the assumptions used here are the same as those used in that paper.
The reference to that paper are proceeded with a M, so that (M:$1$) will refer to the equation
$(1)$ in \cite{Ale20}.

\section{Special functions used in the paper}\label{app:Spec_func}
%
\subsection*{Modification of the harmonic numbers and Lerch transcendent function}
For any integer number $m$, we denote by $\hN{}(m)$ the ${m}$-th harmonic number. We recall here that 
\begin{equation}\label{eq:def_propHarm}
\hN{}(m) = \psi(m+1)+\gamma = \sum_{j=1}^m \frac{1}{j} = \sum_{j=0}^{m-1} \int_0^1 t^j\,dt = 
\int_0^1 \frac{1-t^m}{1-t} \,dt,  
\end{equation}
where $\psi$ is the derivative of the logarithm of gamma function (also called \emph{digamma} function). The constant $\gamma $
is the Euler–Mascheroni constant throughout the whole paper.
The function $\hN{}$ can be extended
therefore to the real non-negative numbers, by setting $\hN{1}(x) = \int_0^1 \frac{1-t^x}{1-t} \,dt$,
which is known as the integral representation given by Euler.

\begin{definition}\label{defn:p-mod}
For $0 \leq p \leq 1$, $x\geq 0$, we define the \emph{$p$-modification of the harmonic numbers} 
$\hN{p}(x)$, where 
\[
\hN{p}(x) = \int_0^1 \frac{1-(1- p + pt)^x}{1-t} \,dt , 
\]
\end{definition}


The function $\hN{p}(x)$ has the following properties
\begin{itemize}
\item $\hN{p}(0)=0$, $\hN{0}(x)=0$, $\hN{p}(1)=p$ and $\hN{1}(x)=\hN{}(x)$ by definition;
\item with two changes of integration variable $z = (1-p(1-t))$ and $z = (1-p)e^{-w}$, we 
we may rewrite $\hN{p}(x)$ 
as
\begin{equation}\label{eq:PHarm3}
\begin{aligned}
\hN{p}(x)  & = \int_{1-p}^1 \frac{1-z^x}{1-z} \,dz = \psi(x+1) +\gamma -\int_0^{1-p} \frac{1-z^x}{1-z} \,dz
\\
& = \psi(x+1) +\gamma +\log p +\int_0^{1-p} \frac{z^x}{1-z} \,dz
\\
& = \psi(x+1) +\gamma +\log p +(1-p)^{x+1}\int_0^{\infty} \frac{e^{-w(x+1)}}{1-(1-p)e^{-w}} \,dw
\\
& = \psi(x+1) +\gamma +\log p +(1-p)^{x+1}\, \Phi(1-p,1,x+1) 
,
\end{aligned}
\end{equation}
where $\Phi$ is the \emph{Lerch transcendent function}, see \cite{OLBC10}, and the last equality is a consequence of
the following equation, valid for $m\in\mathbb{N}$ and $z=(1-p)$:
\[
\Phi\left(z,s,a\right)=z^{m}\Phi\left(z,s,a+m\right)+\sum_{n=0}^{m-1}\frac{z^{%
n}}{(a+n)^{s}}.
\]
\item By \eqref{eq:PHarm3}, $\hN{p}(x)$ is strictly increasing and continuous, both as a function of $x$ and $p$.
In addition, for any $p>0$, $\lim_{x\to\infty}\hN{p}(x)= +\infty$, and hence $\hN{p}:[0,+\infty)\to[0,+\infty)$ is
an isomorphism (continuous invertible function, with continuous inverse function).
Its inverse function $(\hN{p})^{-1}:[0,+\infty)\to[0,+\infty)$ is hence well-defined and it is used in the paper.
\end{itemize}
The Lerch transcendent function appears also in the derivatives of $\hN{p}$. Denote by 
\[
\Phi_{1} = \Phi(1-p,1,x+1), \qquad 
\Phi_{2} = \Phi(1-p,2,x+1), \qquad 
\Phi_{3} = \Phi(1-p,3,x+1), 
\]
and note that $\Phi_{n+1} = - n \partial \frac{\Phi_n}{\partial x}  $;
by \eqref{eq:PHarm3} we get
\begin{equation}\label{eq:PHarmDer}
\begin{aligned}
\hN{p}'(x)  & = 
\partial \frac{\psi(x+1) +\gamma +\log p +(1-p)^{x+1}\cdot \Phi(1-p,1,x+1)}{\partial x}  
\\
& = 
\psi_1(x+1) + (1-p)^{x+1}(\log(1-p) \cdot \Phi_{1} - \Phi_{2} )
\\
\hN{p}''(x)  
& = 
\psi_2(x+1) + (1-p)^{x+1}((\log(1-p))^2 \cdot \Phi_{1} -2 \log(1-p) \cdot \Phi_{2} + 2 \Phi_{3} ).
\end{aligned}
\end{equation}

\subsection*{Product representation and incomplete Gamma function}
For what concerns the infinite product representation of the Gamma function
\[
\Gamma(z) = \lim_{K \to \infty}
\frac{e^{-\gamma z}}{z} \prod_{k=1}^K \left(1 + \frac{z}{k}\right)^{-1} e^\frac{z}{k},
\qquad z \neq -1, -2, \ldots,
\]
given by Schl\"{o}milch in 1844 and Newman in 1848, if we evaluate it in $z = \pm t$,
we obtain
\begin{equation}\label{eq:Gamma}
\Gamma(1-t) e^{-\gamma t} =
\prod_{j=1}^\infty \frac{e^{-\frac{t}{j}}}{1 - \frac{t}{j}} ,
\ t\in (0,1) ,
\qquad
\Gamma(1+t) e^{\gamma t} = 
\prod_{j=1}^\infty \frac{e^{\frac{t}{j}}}{1 + \frac{t}{j}} ,
\ t> 0 .
\end{equation}


\bigskip

Finally, for $x>0$, we denote by $E_1(x)$ the \emph{exponential integral} (or \emph{incomplete gamma function}).
As shown in \cite[p.~229, 5.1.20]{S:AbsSteg64}, we have that
\begin{equation}\label{eq:E1Def}
E_1(x) = \int_x^\infty \frac{e^{-t}}{t} \,dt <
e^{-x}\,\ln\!\Big( 1+\frac{1}{x} \Big) .
\end{equation}
Note that, if $p\in(0,1)$ and $t = -\ln(1-p)w$,
\[
E_1(x) = \int_x^\infty \frac{e^{-t}}{t} \,dt = \int_{-\tfrac{x}{\ln(1-p)}}^\infty \frac{(1-p)^w}{w} \, dw .
\]
We will make use of the very well known formula $-\ln(p) = \sum_{j=1}^\infty \frac{(1-p)^j}{j} $.
To bound the tail of the series, we immediately obtain
by \eqref{eq:E1Def} that, for any $x>0$,
\begin{multline}\label{eq:PHarm4.2}
\sum_{j=0}^\infty \frac{(1-p)^{x+j+1}}{x+j+1} 
\leq  \int_{x}^\infty \frac{(1-p)^w}{w} \, dw 
= E_1(-x\ln({1-p})) 
\\
< 
e^{x\ln({1-p})}\,\ln\!\Big( 1-\frac{1}{x\ln({1-p})} \Big).
\end{multline}
The next representation lemma is used both in the analytical and in the numerical part of the paper.
\begin{lemma}\label{lem:maxPbls}
Let $x>0$ be fixed. Then the functions
\begin{align*}
g_+(t) & = (x+\gamma) t - \ln \Gamma(1-t)  , \qquad t \in (0,1)
\\
g_-(t) & = (x-\gamma) t - \ln \Gamma(1+t)  , \qquad t > 0
\end{align*}
attain their (strictly positive) maxima at the points $t_+ = 1 - \psi^{-1}(-x-\gamma)$ and
$t_- = \psi^{-1}(x-\gamma) - 1$, respectively.
\end{lemma}
\begin{proof}
We give the proof for $g_+$, since the same arguments apply to $g_-$. We have
\begin{itemize}
\item $g_+(t)$ is concave, since $\ln\Gamma(1-t)$ is a convex analytic function on $(0,1)$;
\item $g_+(0) = \ln\Gamma(1)=0 $, $g_+'(0) = (x+\gamma) + \psi(1) = x >0$;
\item $\lim_{t\to 1} g_+(t) = -\infty$;
\end{itemize}
and hence 
the maximum of $g_+$ on $(0,1)$ is strictly positive. 
The maximum point $t_+$ is attained when $g_+'(t_+)=0$, that is
when $(x+\gamma)+\psi(1+t_+) = 0$.
\end{proof}

\section{Proof of some technical results of \cite{Ale20}}

\subsection{Proof of $0< {Y}_{r\,c}  - \bar{Y}_{r\,c} < \tfrac{2^{-z_0}}{\lambda_0}$ in \cite[Lemma~\ref{cor:Y_rc}]{Ale20}}\label{sec:Ybounds}
We recall here that
\[
2^{-z_0} Z(c,o) = \sum_{z = 1}^{z_0} s^*_{z} 2^{-z} , \qquad
\bar{Z}(o,c) = \sum_{z = 1}^{\infty} s^*_{z} 2^{-z}
\]
and with 
\[
{Y}(o,c) = X(c,o) - \log_2 (1+ 2^{-z_0} Z(c,o) ), \qquad
\bar{Y}(o,c) = \big( X(c,o) - \log_2 (1+ \bar{Z}(o,c) ),
\]
by definition of ${Y}_{r\,c} $ and $\bar{Y}_{r\,c}$, we get
\begin{align*}
{Y}_{r\,c}  - \bar{Y}_{r\,c} & = 
\max_{\{o \colon R(o,c)= r\}} ( {Y}(o,c) - \bar{Y}(o,c) ) 
\\
& = 
\max_{\{o \colon R(o,c)= r\}} \big( - \log_2 (1+ 2^{-z_0} Z(c,o) ) + 
\log_2 (1+ \bar{Z}(c,o) ) 
\big) 
\\
& = 
\max_{\{o \colon R(o,c)= r\}} 
\log_2 \Big( 1 + \frac{2^{-z_0}\sum_{z=z_0+1}^{\infty} s^*_{z} 2^{-z}}{
1 + \sum_{z=1}^{z_0} s^*_{z} 2^{-z} } \Big).
\end{align*}
Now, note that 
\[
0 < \frac{\sum_{z=1}^{\infty} s^*_{z} 2^{-z}}{
1 + \sum_{z=1}^{z_0} s^*_{z} 2^{-z} } < 1
\]
and then, since $\ln(1+x)<x$ for $x>0$,
\[
0 < \log_2 \Big( 1 + \frac{2^{-z_0}\sum_{z=z_0+1}^{\infty} s^*_{z} 2^{-z}}{
1 + \sum_{z=1}^{z_0} s^*_{z} 2^{-z} } \Big) <
\frac{2^{-z_0}}{\ln(2)} = \frac{2^{-z_0}}{\lambda_0}. 
\]

\subsection{Detailed proof of  \cite[Theorem~\ref{thm:ris1}]{Ale20}}\label{sec:SMproof}

\noindent \textbf{First step}.
By \cite[Lemma~\ref{cor:Y_rc}]{Ale20}, it is possible to calculate the moment-generating function of $\bar{\mathcal{Y}} $, conditioned
on $\{\boldsymbol{m}_c , c= 1,\ldots,c_0\}$.
In fact, since 
\begin{equation}\label{eq:barYrcDef}
\bar{Y}_{r\,c} = 
\max_{\substack{o_1, \ldots,o_{m_{r\,c}}\colon R(c,o_j)=r\\ o_j \text{ different objects}}} \big( \bar{Y}(o_j,c) \big) ,
\end{equation}
it is well known \cite{TbCited} that the moment generating function of the max of exponential random variables is
\[
E(e^{s\bar{Y}_{r\,c}} | \{\boldsymbol{m}_c , c= 1,\ldots,c_0\} ) = \prod_{j=1}^{m_{r\,c}} (1 - \tfrac{s}{\lambda_0 j} )^{-1},
\qquad 0<s<\lambda_0
\]
which implies, for $0<s<{c_0}2^{r_0} \lambda_0$,
\begin{multline*}
E(e^{s\bar{\mathcal{Y}}} | \{\boldsymbol{m}_c , c= 1,\ldots,c_0\}) = 
E \Big(e^{s
\sum_{c=1}^{{c_0}}
\sum_{r=1}^{2^{r_0}}
\tfrac{\bar{Y}_{r\,c}}{{c_02^{r_0}}}} \Big| \{\boldsymbol{m}_c , c= 1,\ldots,c_0\}\Big) 
\\
= 
\prod_{c=1}^{{c_0}}
\prod_{r=1}^{2^{r_0}}
\prod_{j=1}^{m_{r\,c}} (1 - \tfrac{s}{j {c_0}2^{r_0} \lambda_0} )^{-1} .
\end{multline*}
Again, by \eqref{eq:barYrcDef}
\begin{equation}\label{eq:barYCond}
E ( Y_{r\,c} | \{\boldsymbol{m}_c , c= 1,\ldots,c_0\} ) = 
\sum_{j=1}^{m_{r\,c}} \frac{1}{j \lambda_0}
\end{equation}
which means that 
\begin{equation*}
E({\bar{\mathcal{Y}}} | \{\boldsymbol{m}_c , c= 1,\ldots,c_0\}) =
\frac{1}{c_02^{r_0}\lambda_0}
\sum_{c=1}^{{c_0}}
\sum_{r=1}^{2^{r_0}}
\sum_{j=1}^{m_{r\,c}} \frac{1}{j}.
\end{equation*}
Then, conditioned on $\{\boldsymbol{m}_c , c= 1,\ldots,c_0\}$, 
the Chernoff bound for the first inequality that concerns $\alpha_+$ may be computed as
\begin{multline*}
P\Big( \bar{\mathcal{Y}} \geq E(\bar{\mathcal{Y}}|\{\boldsymbol{m}_c , c= 1,\ldots,c_0\})+\frac{h_d}{\lambda_0} \Big|
\{\boldsymbol{m}_c , c= 1,\ldots,c_0\}
\Big) \\
\begin{aligned}
& \leq \min_{s>0}
e^{-s\big(
{
 E(\bar{\mathcal{Y}}|\{\boldsymbol{m}_c , c= 1,\ldots,c_0\})+\frac{h_d}{\lambda_0} 
 } \big)
 }   
E(e^{s\bar{\mathcal{Y}}} | \{\boldsymbol{m}_c , c= 1,\ldots,c_0\}) 
\\ &
= \min_{s>0}
e^{-\tfrac{{h_d}}{\lambda_0}s}
\prod_{c=1}^{{c_0}}
\prod_{r=1}^{2^{r_0}}
\prod_{j=1}^{m_{r\,c}} 
\frac{
e^{-\tfrac{s}{j
c_02^{r_0}\lambda_0}}
}{
1 - \tfrac{s}{j {c_0}2^{r_0} \lambda_0}
}.
\end{aligned}
\end{multline*}
Define $t = \frac{s}{{c_0}2^{r_0} \lambda_0}$. 
Since $\frac{\exp^{-t}}{1-t}\geq 1$ for any $t<1$, then for $t \in
(0, 1)$, the abobe relation continues as
\begin{align*}
\min_{s>0}
e^{-\tfrac{{h_d}}{\lambda_0}s}
\prod_{c=1}^{{c_0}}
\prod_{r=1}^{2^{r_0}}
\prod_{j=1}^{m_{r\,c}} 
\frac{
e^{-\tfrac{s}{j
c_02^{r_0}\lambda_0}}
}{
1 - \tfrac{s}{j {c_0}2^{r_0} \lambda_0}
}
& = 
\min_{t>0}
e^{-t{h_d} {c_0}2^{r_0}}
\prod_{c=1}^{{c_0}}
\prod_{r=1}^{2^{r_0}}
\prod_{j=1}^{m_{r\,c}} 
\frac{
e^{-\tfrac{t}{j}} } {
1 - \tfrac{t}{j}
}
\\
& 
\leq
\min_{t>0}
e^{-t{h_d} {c_0}2^{r_0}}
\prod_{c=1}^{{c_0}}
\prod_{r=1}^{2^{r_0}}
\prod_{j=1}^{\infty} 
\frac{
e^{-\tfrac{t}{j}} } {
1 - \tfrac{t}{j}
}
\\
& 
=
\min_{t>0}
e^{-t{h_d} {c_0}2^{r_0}}
\prod_{c=1}^{{c_0}}
\prod_{r=1}^{2^{r_0}}
\Big(\Gamma(1-t) e^{-\gamma t} \Big)
\\
& 
= \exp\Big( -{c_0}2^{r_0} \max_{t\in(0,1)}  \big[ (h_d+\gamma) t - \ln \Gamma(1-t) \big]\Big).
\end{align*}
The relevant aspect of the last expression is that it does not depend on $\{\boldsymbol{m}_c , c= 1,\ldots,c_0\}$, and
what remains to prove is that the maximum of 
\( (h_d+\gamma) t - \ln \Gamma(1-t) \) on $(0,1)$
is attained at
\(
t_+ = 1-\psi^{-1}(-h_d-\gamma)
\).
This is obvious, since $- \ln \Gamma(1-t)$ is a concave function with derivative in zero equal to
$-\gamma$, its limit is $-\infty $ as it approches $1^-$ and the digamma function $\psi$ is the derivative of
the logarithm of the gamma function.

\bigskip 

The proof of the second inequality that concerns $\alpha_-$ may be done with the same ideas.
In fact, the Chernoff bound may be uniformly bounded by
\begin{multline*}
P\Big( \bar{\mathcal{Y}} \leq E(\bar{\mathcal{Y}}|\{\boldsymbol{m}_c , c= 1,\ldots,c_0\})-\frac{h_u}{\lambda_0} \Big|
\{\boldsymbol{m}_c , c= 1,\ldots,c_0\} \Big)
\\
\leq \min_{t>0}
\exp\Big( -{c_0}2^{r_0} \max_{t>0}  \big[ (h_u-\gamma) t - \ln \Gamma(1+t) \big]\Big).
\end{multline*}
Now, it is sufficient to note that $- \ln \Gamma(1+t)$ is a concave function with derivative in zero equal to
$\gamma$, its limit is $-\infty $ as it approches $+\infty$ and the digamma function $\psi$ is again the derivative of
the logarithm of the gamma function.

\bigskip 

\noindent \textbf{Second step}. Starting from \eqref{eq:barYCond}, note that
the $m_{r\,c}$-armonic number may be represented in the following way:
\[
\sum_{j=1}^{m_{r\,c}} \frac{1}{j} = 
\int_0^1 \frac{1-v^{m_{r\,c}}}{1-v} \,dv = \hN{1}(m_{r\,c}).
\]
Recall that, by \cite[Lemma~\ref{cor:Y_rc}]{Ale20}, $m_{r\, c}$
is distributed as a binomial distribution, with $F_0$ trials and probability $p_0={2^{-{r}_0}}$.
Then
\begin{align*}
\lambda_0 E(Y_{r\, c} )  
& = 
\lambda_0 E( E(Y_{r\, c}| \{\boldsymbol{m}_c , c= 1,\ldots,c_0\} ) )  = E( \hN{1}(m_{r\, c}) )
\\
& 
= \sum_{m=0}^{F_0} \hN{1}(m) \binom{{F_0}}{m} {p_0}^m (1-{p_0})^{{F_0}-m} \\
& = \sum_{m=0}^{F_0}  \Big(\int_0^1 \frac{1-v^m}{1-v} \,dv\Big)  \binom{{F_0}}{m} {p_0}^m (1-{p_0})^{{F_0}-m} \\
& = \int_0^1 \frac{1}{1-v} \Big(\sum_{m=0}^{F_0} (1-v^m)\binom{{F_0}}{m} {p_0}^m (1-{p_0})^{{F_0}-m}  \Big) dv  \\
& = \int_0^1 \frac{1}{1-v} \Big(\sum_{m=0}^{F_0} \binom{{F_0}}{m} {p_0}^m (1-{p_0})^{{F_0}-m}  \\
& \qquad\qquad - \sum_{m=0}^{F_0} \binom{{F_0}}{m} ({p_0}v)^m (1-{p_0})^{{F_0}-m}  \Big) dv  \\
& = \int_0^1 \frac{1- (1-{p_0} + {p_0}v)^{F_0} }{1-v} \,dv  = \hN{p_0} (F_0) .
\end{align*}
Then, by linearity, we conclude that
$E(\bar{\mathcal{Y}}) = \frac{\hN{p_0}(F_0)}{\lambda_0}$.

\section{Lower and upper bounds of some numerical problems} \label{app:boundScomp}
\subsection{Bounds of $y = \psi(x)$}\label{sec:boundPsi}
As shown in \cite[Example 2.1]{S:DS2016}, we may bound $\psi$ from below in the following way.
The Jensen inequality for $U\sim U({x-\tfrac{1}{2}},{x+\tfrac{1}{2}})$ shows that, for $x>\tfrac{1}{2}$,
\[
\frac{1}{x} 
= \frac{1}{E[U]} < E \Big[\frac{1}{U} \Big] = \int_{x-\tfrac{1}{2}}^{x+\tfrac{1}{2}} \frac{1}{t} \,dt
= \ln (x+\tfrac{1}{2})  - \ln (x-\tfrac{1}{2}) .
\]
By \meqref{eq:psiDef1}, we than have that, for $x>\tfrac{1}{2}$,
\[
\psi(x) - \ln (x-\tfrac{1}{2}) > \psi(x+1) - \ln (x+\tfrac{1}{2}) > \cdots > \mathop{\lim\inf}_{t\to\infty} (\psi(t) - \ln (t-\tfrac{1}{2} )),
\]
and since $\psi(t) = \log(t)+o(1) = \log(t-\tfrac{1}{2})+o(1)$, the last expression is zero, and hence
\[
y = \psi(x) > \ln (x-\tfrac{1}{2}) , \qquad \text{for any }x>\tfrac{1}{2}.
\]
With the same spirit of this example, since
\[
\ln (x+1)  - \ln (x) = \int_{x}^{x+1} \frac{1}{t} \,dt< \frac{1}{x} , \qquad
\forall x> 0,
\]
we obtain that
\[
\psi(x) - \ln (x) < \psi(x+1) - \ln (x+1) < \cdots < \mathop{\lim\sup}_{t\to\infty} (\psi(t) - \ln (t)) =0,
\]
and hence, we may state that
\begin{equation*}
\ln (x-\tfrac{1}{2}) < y < \ln (x) , \qquad e^{y} < x < e^{y} + \tfrac{1}{2}, 
\qquad \forall x > \tfrac{1}{2}, \forall y=\phi(x).
\end{equation*}

\subsection{Bounds of $y = \hN{p}(x)$}\label{sec:boundHnP}
For what concerns the bounds for $\hN{p}$, by \eqref{eq:PHarm3}, we immediately get
\[
\psi(x+1)+\gamma+ \ln p  \leq \hN{p}(x) \leq  
\psi(x+1)+\gamma,
\]
and hence, by \meqref{eq:psiAppr},
\begin{equation}\label{eq:center_eq1}
\frac{\exp ( \hN{p}(x) - \gamma )}{p} - \frac{1}{2} 
\geq x \geq
\exp ( \hN{p}(x) - \gamma ) - 1 .
\end{equation}
A better estimation 
for the lower bound 
can be found for  $x > -\tfrac{1}{(e-1)\ln(1-p)}$. 
To simplify the notations, set $d_0 = -{\ln(1-p)}$, so that the assumption
$x > -\tfrac{1}{(e-1)\ln(1-p)}$ becomes the more readable 
$x d_0 > \tfrac{1}{e-1}$. 
We are going to show that, under this hypothesis, we have
\begin{equation}\label{eq:center_eq}
\frac{A}{p} - \frac{1}{2} 
\geq
x \geq
\begin{cases}
\frac{A}{p} - e +\tfrac{1}{\ln(1-p)} & 
\text{if }A > p(\tfrac{1}{2} -\tfrac{1}{(e-1)\ln(1-p)}) ; \\
A - 1 & 
\text{otherwise};
\end{cases}
\end{equation}
where $A = \exp ( \hN{p}(x) - \gamma )$.
To prove \eqref{eq:center_eq}, we use the relation $\frac{1}{1-z} = \sum_{j=0}^\infty z^j$,
valid for $|z|<1$, in \eqref{eq:PHarm3}. We obtain
\begin{equation*}
\begin{aligned}
\hN{p}(x) 
& = \psi(x+1) +\gamma +\log p +\int_0^{1-p} \frac{z^x}{1-z} \,dz
\\
& = \psi(x+1) +\gamma +\log p + \int_0^{1-p} \sum_{j=0}^\infty z^{x+j} \,dz
\\
& = \psi(x+1) +\gamma +\log p + \sum_{j=0}^\infty \frac{(1-p)^{x+j+1}}{x+j+1} \,dz
,
\end{aligned}
\end{equation*}
which can be combined with \eqref{eq:PHarm4.2}, yielding
\begin{multline}\label{eq:PHarm4.3}
\hN{p}(x) - (\psi(x+1) +\gamma +\log p )
  < 
e^{x\ln({1-p})}\,\ln\!\Big( 1-\frac{1}{x\ln({1-p})} \Big) \\
< e^{x\ln({1-p})} 
\leq \frac{1}{1 - x\ln({1-p})},
\end{multline}
where the last inequality is a consequence of the fact that $\exp(x) \leq \tfrac{1}{1-x}$ for $x<1$.

Now, we define the positive quantity
$d_1 = e-1+\tfrac{1}{d_0}>0$ and we note that the function $g:[\frac{1}{d_0(e-1)},\infty)\to \mathbb{R}$
so defined
\[
g(x) = \frac{d_1}{d_1+1+x}
- \frac{1}{1 + xd_0}  = \frac{x(d_0d_1-1)-1}{(d_1+1+x)(1 + xd_0)} 
\]
is strictly positive whenever $x(d_0d_1-1)-1> 0$, 
or, in other terms, when $d_1 > \tfrac{1 + x}{d_0 x}$. We now prove that 
this fact implies that $g(x)>0$ under our assumption $x > \tfrac{1}{d_0(e-1)} $.

In fact, since $\tfrac{1 + y}{d_0 y}$
is decreasing in $y>0$, then, as $x >\tfrac{1}{d_0(e-1)}$
we have
\[
x > \tfrac{1}{d_0(e-1)} \quad \Longrightarrow \quad 
d_1 = \tfrac{d_0(e-1)+1}{d_0} = \frac{1 + \tfrac{1}{d_0(e-1)} }{d_0 \tfrac{1}{d_0(e-1)}} > 
\tfrac{1 + x}{d_0 x} \quad \Longrightarrow \quad
g(x) >0,
\]
or, in other terms,
\[
x > \tfrac{1}{d_0(e-1)} \quad \Longrightarrow \quad 
\frac{d_1}{d_1+1+x}
> \frac{1}{1 + xd_0} = \frac{1}{1 - x\log(1-p)} .
\]
Since $\frac{x}{1+x}< \ln(1+x)$ for $x>0$,
we then have that, when $x > \tfrac{1}{d_0(e-1)}$,
\begin{multline}\label{eq:PHarm4.4}
\frac{1}{1 - x\ln({1-p})} 
<
\frac{d_1}{d_1+1+x}  
=
\frac{\frac{d_1}{x+1}}{1+\frac{d_1}{x+1}} 
<
\log \Big(1 + \frac{d_1}{x+1} \Big) 
\\
= \ln \Big( \frac{x + 1 + d_1}{x+1} \Big) 
= \ln ({x + e -\tfrac{1}{\ln(1-p)}}) 
 - \ln ({x+1}) .
\end{multline}
By combining together \eqref{eq:PHarm4.3} and \eqref{eq:PHarm4.4}
we obtain
\[
\hN{p}(x) - (\psi(x+1) +\gamma +\log p ) < 
\ln ({x + e -\tfrac{1}{\ln(1-p)}}) 
 - \ln ({x+1}) ,
\]
that together with \meqref{eq:psiAppr} yields 
\begin{align*}
\hN{p}(x) - \gamma -\log p & < 
\psi(x+1) - \ln ({x+1}) + \ln ({x + e -\tfrac{1}{\ln(1-p)}}) \\
& < 
\ln ({x + e -\tfrac{1}{\ln(1-p)}}) .
\end{align*}
Set $A = \exp ( \hN{p}(x) - \gamma )$. The above inequality, exponentiated, gives
\[
\frac{A}{p} - e +\tfrac{1}{\ln(1-p)} < x,
\]
that, again by \eqref{eq:center_eq1}, is valid at least when
\[
x > -\tfrac{1}{(e-1)\ln(1-p)} \quad \Longrightarrow \quad  A > p( x + \tfrac{1}{2}) > 
p(\tfrac{1}{2} -\tfrac{1}{(e-1)\ln(1-p)}) .
\]

\subsection{Bounds of $y= (x-\gamma) t(x) - \ln \Gamma(1+t(x)) $, where $t(x)= \psi^{-1}(x-\gamma)-1$}\label{sec:boundAm}
For what concerns the bounds in this problem, we start by recalling that,
as shown in \cite{LN13} (see also \cite[Equation (3.112)]{Qi2010}), for any $t > 0$, 
\begin{equation}\label{eq:Qi}
-\gamma t < \ln \Gamma (1 + t) < t \psi(t+1) .
\end{equation}
When this chain of inequalities is evaluated in $t= t(x)$, we obtain 
\begin{align}
-\gamma t(x) - \ln \Gamma (1 + t(x)) <0
& \quad \Longrightarrow \quad
y < x t(x) \label{eq:imediate}
\\
\ln \Gamma (1 + t(x)) < t \psi(\psi^{-1}(x-\gamma)-1+1)
& \quad \Longrightarrow \quad
y >0. \notag
\end{align}
The upper bounds for $x$ may be found in the following way. We recall that Lemma~\ref{lem:maxPbls}
states that
\[
y = \max_{t>0}\big[ (x-\gamma) t - \ln \Gamma(1+t) \big]. 
\]
Then, by \eqref{eq:Qi}, 
\begin{equation}\label{eq:disUB+}
y > \max_{t>0}\big[ (x-\gamma - \psi(t+1) ) t \big] .
\end{equation}
The second expression may be evaluated in \(
t_0 = \psi^{-1}( -\gamma + \tfrac{x}{2}) -1
\), so that we get
\begin{align}
y & > 
(x-\gamma - \psi(t_0+1) ) t_0 \notag
\\
& 
= \frac{x}{2} \big( \psi^{-1}( -\gamma + \tfrac{x}{2}) -1 \big) \label{eq:eval_halph}
\\
& 
= \frac{x}{2} \big( \psi^{-1}( -\gamma + \tfrac{x}{2}) -\psi^{-1}( -\gamma) \big). \notag
\end{align}
The Mean Value Theorem ensures the existence of $x_0 \in (0, \tfrac{x}{2})$ such that
\[ 
\psi^{-1}( -\gamma + \tfrac{x}{2}) -\psi^{-1}( -\gamma)  = \frac{x}{2} {d\frac{\psi^{-1}(-\gamma + t)}{dt}\Big|_{t=x_0}} ,
\]
and by the the formula of the derivative of the inverse function, since the Trigamma function $\psi_1(t) = d\frac{\psi(t)}{dt}$ is a decreasing function
with $\psi_1(1) = \frac{\pi^2}{6}$,
\[
 {d\frac{\psi^{-1}(-\gamma + t)}{dt}\Big|_{t=x_0}} = 
 \frac{1}{\psi_1(\psi^{-1}(-\gamma + x_0))}
>   \frac{1}{\psi_1(\psi^{-1}(-\gamma))}
=   \frac{1}{\psi_1(1)}
=  \frac{1}{\frac{\pi^2}{6}}.
\]
Summing up,
\begin{equation}\label{eq:UB+1}
y > \frac{x}{2} \Big( \frac{x}{2} \frac{1}{\frac{\pi^2}{6}} \Big) = \frac{3}{2}\frac{x^2}{\pi^2} 
\quad \Longrightarrow \quad
x < \pi\sqrt{\frac{2}{3}  y}.
\end{equation} 
For $x\geq \frac{3}{2}$, which is always true if $y\geq \frac{3}{2}\cdot t(\frac{3}{2}) =3$ by \eqref{eq:imediate}, 
a better estimates may be found 
if we bound the second part of \eqref{eq:eval_halph}. In fact, since $\frac{x}{2}\geq \frac{3}{4}$,
by \meqref{eq:psiAppr} we obtain 
\[
y > 
\frac{3}{4} \big( \psi^{-1}( -\gamma + \tfrac{x}{2}) -1 \big) 
> \frac{3}{4} \big( e^{-\gamma + \tfrac{x}{2}} -1 \big) 
\]
%
%
which completes the upper bound for $x$ given in \eqref{eq:UB+1}, obtaining
\begin{equation}\label{eq:-bounds}
x < 
\begin{cases}
\pi\sqrt{\frac{2}{3}  y} , & \text{if $y <3$;}
\\
2 (\log (\frac{4}{3}y+1) + \gamma) 
, & \text{if $y \geq 3$.}
\end{cases}
\end{equation}
The upper bounds for $x$ may be found with similar ideas in both the cases $y\geq 3$ and $y<3$.
By \eqref{eq:imediate}, the Mean Value Theorem ensures the existence of 
$x_0 \in (0, x)$ such that, when $y<3$ 
\[
0 < y  < x t(x) = x(\psi^{-1}(x-\gamma)-1) = x^2 \frac{1}{\psi_1(\psi^{-1}(x_0-\gamma))}
<  x^2 \frac{1}{\psi_1(\psi^{-1}(\pi\sqrt{2}-\gamma))},
\]
the last inequality being a consequence of \eqref{eq:-bounds}, since, for $y<3$, we have $x \leq \pi\sqrt{2}$. 
For $y\geq 3$, starting from \eqref{eq:imediate},
by \eqref{eq:psiAppr}, we obtain
\[
0 < y  < x t(x) = x(\psi^{-1}(x-\gamma)-1) < x \Big( \exp(x-\gamma) - \frac{1}{2}\Big) < \Big( \exp(\tfrac{3}{2}x-\gamma) - \frac{1}{2}\Big) ,
\]
which gives the lower bound for $x$ in \meqref{eq:boundsAmAll} for $y\geq 3$.

\subsection{Bounds of $y= (x+\gamma) t(x) - \ln \Gamma(1-t(x)) $, where $t(x)= 1-\psi^{-1}(-x-\gamma)$}\label{sec:boundAp}
The inversion formula for the Gamma function, valid for $t\in(0,1)$, gives
\[
\Gamma(1-t) \Gamma(t) t= {\pi \over \sin{(\pi t)}} t 
\qquad \Longleftrightarrow \qquad 
 \ln\Gamma(1-t)  = \ln\Big({\pi t \over \sin{(\pi t)}}\Big) -  \ln\Gamma(1+t),  
\]
that, together with \eqref{eq:Qi}, 
yealds
\begin{equation}\label{eq:inverGammaQi}
-t \psi(t+1) +\ln\Big({\pi t \over \sin{(\pi t)}}\Big)  < \ln\Gamma(1-t)  < \ln\Big({\pi t \over \sin{(\pi t)}}\Big) + \gamma t .
\end{equation}
We recall that Lemma~\ref{lem:maxPbls}
states that
\[
y = \max_{t\in(0,1)}\big[ (x+\gamma) t - \ln \Gamma(1+t) \big],
\]
that, combined with the right-hand inequality of \eqref{eq:inverGammaQi} gives
\[
y 
> \max_{t\in(0,1)}\big[ x t + \ln\Big({\sin(\pi t)\over \pi t} \Big)\big] .
\]
Since $\ln(y) > \tfrac{y-1}{y}$ and  (see \cite{AZ10}), 
\[
\frac{\pi}{\sin(\pi t)} = \frac{1}{t} + \sum_{n=1}^\infty \frac{(-1)^n\,2t}{t^2-n^2},
\]
then
\[
y > \max_{t\in(0,1)}\big( x t -  \frac{2t^2}{1-t^2} \big) .
\]
Let $t_0= t_0(x) \in (0,1)$ be defined in the following way:
\[
\frac{x}{2} =  \frac{2t_0}{1-t_0^2}  \qquad \Longleftrightarrow \qquad t_0 = 2\frac{ \sqrt{\big(\frac{x}{2}\big)^2 + 1} - 1}{x} ,
\]
then 
\[
y >  x t_0 -  t_0\frac{2t_0}{1-t_0^2}  = x \frac{t_0}{2} = \sqrt{\Big(\frac{x}{2}\Big)^2 + 1} - 1,
\]
and hence
\begin{equation}\label{eq:+ub1}
x < 2 \sqrt{(y+1)^2-1}.
\end{equation}
For what concerns the lower bound for $x$, 
if we take into account the reflection formula for the digamma function
\[
\psi(1-t)-\psi(t)=\pi \cot \pi t 
\qquad \Longrightarrow \qquad 
\psi(1 +t) = \psi(t) +\frac{1}{t} = \psi(1-t) - \pi \cot (\pi t) + \frac{1}{t}
\]
together with 
the left inequality in \eqref{eq:inverGammaQi}, we obtain
\begin{align*}
\ln\Gamma(1-t) & > -t \psi(t+1) +\ln\Big({\pi t \over \sin{(\pi t)}}\Big)  
\\
& = -t \Big( \psi(1-t) - \pi \cot (\pi t) + \frac{1}{t}\Big) +\ln\Big({\pi t \over \sin{(\pi t)}}\Big) .
\end{align*}
We will make use of this inequality, motivated by the fact that our problem is
\[
y =  (x+\gamma) t(x) - \ln \Gamma(1-t(x)) , \qquad 
-\psi(1 -t(x)) = (x+\gamma),
\]
which implies
\begin{align} \notag
y  & =  (x+\gamma)  t(x) - \ln \Gamma(1-t(x)) 
\\ \notag
& =  - \psi(1 -t(x))  t(x) - \ln \Gamma(1-t(x)) 
\\
& 
< 1 - \pi t(x) \cot (\pi t(x))  + \ln\Big({\sin{(\pi t(x))}  \over \pi t(x) }\Big) . \label{eq:+lb0}
\end{align}
Now, for $t \in (0,1)$, the following identities hold
\[
{\sin{(\pi t)}  \over \pi t } = \prod_1^{\infty} \Big( 1 - \frac{t^2}{ n^2} \Big)\,,
\qquad
\pi \cdot \cot(\pi t) = \frac{1}{t} + \sum_{n=1}^\infty \frac{2t}{t^2-n^2}\,, 
\]
(see \cite{AZ10}). The first identy may be used to bound the last term in \eqref{eq:+lb0}:
\begin{multline*}
\ln\Big({\sin{(\pi t)}  \over \pi t }\Big) =
\ln(1-t^2) +\sum_{n=2}^\infty \ln\Big( 1 - \frac{t^2}{ n^2} \Big) 
<
\ln(1-t^2) + t^2 -\sum_{n=1}^\infty \frac{t^2}{ n^2} 
\\
=  \ln(1-t^2) + t^2 \Big( 1 -\frac{\pi^2}{6} \Big).
\end{multline*}
For what concerns the term $1 - \pi t \cot (\pi t)$ in \eqref{eq:+lb0}, we obtain
\begin{multline*}
1 - \pi t \cot (\pi t)
= 2t^2 \sum_{n=1}^\infty \frac{1}{n^2-t^2}  
= 2t^2 \Big( \frac{1}{1-t^2} + \sum_{n=2}^\infty \frac{1}{n^2-t^2}  \Big)
\\
< 
2t^2 \Big( \frac{1}{1-t^2} + \sum_{m=1}^\infty \frac{1}{(m+1)^2-1}  \Big)
= 2t^2 \Big( \frac{1}{1-t^2} + \frac{1}{2} \sum_{m=1}^\infty \frac{2}{m(m+2)}\Big)
\\
= 2t^2 \Big( \frac{1}{1-t^2} + \frac{1}{2} \sum_{m=1}^\infty \big(\frac{1}{m} - \frac{1}{m+2}\big)   \Big)
= \frac{2t^2}{1-t^2} + 3t^2 .
\end{multline*}
Combining these two last inequalities in \eqref{eq:+lb0}, since $\log y \leq y-1$, we obtain 
\[
y < \frac{2t(x)^2}{1-t(x)^2} + 3t(x)^2 + \ln(1-t(x)^2) + t(x)^2 \Big( 1 -\frac{\pi^2}{6} \Big) 
< 
\frac{2}{1-t(x)^2} - 2 +t(x)^2 \Big( 3 -\frac{\pi^2}{6} \Big) ,
\]
and hence, if  we define
\[
z = 1-t(x)^2 \in (0,1), \qquad A = \frac{3 -\frac{\pi^2}{6} }{2}\in(0,1), \qquad B = \frac{y }{2} >0
\]
we obtain
\[
Az^2+(B+(1-A))z -1 <0, \qquad z \in(0,1)
\]
which is solved for
\[
0 < z < \frac{-(B+(1-A)) + \sqrt{(B+(1-A))^2 +4A} }{ 2 A}.
\]
Note that, for $B\in(0,\infty)$, the right-hand side of the inequality above belongs to $(0,1)$.
Then, if we define
\[
C = \sqrt{ 1 - \frac{-(B+(1-A)) + \sqrt{(B+(1-A))^2 +4A} }{ 2 A}} \in (0,1),
\]
we have $t(x) = \sqrt{1-z} >C$, or explicitely
\begin{equation}\label{eq:pointFixBound}
1-\psi^{-1}(-x-\gamma) > C.
\end{equation}
Two inequalities on $x$ are consequence of \eqref{eq:pointFixBound} as follows.
By \meqref{eq:psiAppr} we imediately obtain a lower bound 
\[
1- \exp(-(x+\gamma)) >
1-\psi^{-1}(-x-\gamma) > C
\qquad \Longrightarrow \qquad 
x > - \ln (1 - C) - \gamma ,
\]
which is meaningful only for $C\geq 1-\exp(-\gamma)$.
For smaller $C$, we make use of
the Mean Value Theorem, that ensures the existence of $x_0 \in (0, x)$ such that
\[ 
t(x)=\psi^{-1}( -\gamma) - \psi^{-1}( -\gamma - x)  = -x\, {d\frac{\psi^{-1}(-\gamma - t)}{dt}\Big|_{t=x_0}} .
\]
The formula of the derivative of the inverse function gives
\[
- {d\frac{\psi^{-1}(-\gamma - t)}{dt}\Big|_{t=x_0}} = 
 \frac{1}{d\frac{\psi(t)}{dt}\big|_{t = \psi^{-1}(-\gamma - x_0)}}
 <  \frac{1}{d\frac{\psi(t)}{dt}\big|_{t = \psi^{-1}(-\gamma)}}
=  \frac{1}{\frac{\pi^2}{6}},
\]
so that
\[
x > \frac{\pi^2}{6} C .
\]
Summing up
\begin{equation}\label{eq:+bounds}
x > \max\Big( - \ln (1 - C) - \gamma , \frac{\pi^2}{6} C \Big) ,
\end{equation}
that completes \meqref{eq:boundsApAll} 
with the upper bounds for $x$ given in \eqref{eq:+ub1}.

\end{document}